\documentclass[a4paper,11pt,reqno]{article}
\usepackage[T1]{fontenc}
\usepackage{courier}
\usepackage[scaled=0.92]{helvet}
\usepackage{amscd,amsmath,amssymb,amsthm,amsfonts,epsfig,graphics,accents}

\usepackage[percent]{overpic} 				
\usepackage[english]{babel}
\usepackage{amsthm,amsfonts,amssymb,amsmath,subfigure}
\usepackage{hyperref}
\hypersetup{
colorlinks=true,
citecolor=red,
linkbordercolor={1 1 1},
citebordercolor={1 1 1},
}
\usepackage{cleveref}
\usepackage{mathrsfs}
\usepackage{color}
\usepackage{bm}

\usepackage{xcolor} 

\newcommand{\mx}{\mbox}
\newcommand{\rw}{\rightarrow}

\newcommand{\ml}{\mathcal}

\newcommand{\beq}[1]{\begin{equation} \label{#1}}
\newcommand{\eeq}{\end{equation}}
\newcommand{\beqar}{\[ \begin{array}{rcl}}
\newcommand{\eeqar}{\end{array} \]}

\newcommand{\ue}[1]{\boldsymbol{#1}}
\newcommand{\uue}[1]{{\boldsymbol{#1}}}

\providecommand{\ep}{\varepsilon}

\providecommand{\RR}{\mathbb{R}}
\providecommand{\CC}{\mathbb{C}}

\providecommand{\NN}{\mathbb{N}}

\newtheorem{satz}{Theorem}[section]
\newtheorem{prop}{Proposition}[section]

\newtheorem{lem}[satz]{Lemma}

\newtheorem{rem}{Remark}[section]

\usepackage{xcolor}

\newcommand{\ta}{\tilde{\alpha}}
\newcommand{\tb}{\tilde{\beta}}
\newcommand{\thd}{\tilde{h}_2}

\usepackage{xcolor}

\providecommand{\ep}{\varepsilon}

\providecommand{\RR}{\mathbb{R}}
\providecommand{\CC}{\mathbb{C}}

\providecommand{\NN}{\mathbb{N}}

\usepackage[utf8]{inputenc}

\hypersetup{pdfauthor={Fortunati, A. et al.},%
            pdftitle={Title to be defined...},%
            pdfsubject={},%
            pdfkeywords={}%
}

\DeclareMathOperator{\const}{const.}

\DeclareMathOperator{\diag}{diag}

\makeatletter
\renewcommand*{\@fnsymbol}[1]{\ensuremath{\ifcase#1\or *\or \mathsection \or \dagger \or \else \fi}}
\makeatother

\title{\LARGE{\textbf{{Controlled Quasi-Latitudinal Solutions for ultra-fast Spin-Torque Precessional Magnetization Switching}}}}
\date{}
\author{%
Alessandro Fortunati\thanks{E-mail: alessandro.fortunati@unina.it}, Massimiliano d'Aquino\thanks{E-mail: mdaquino@unina.it} $\,$ and Claudio Serpico\thanks{E-mail: serpico@unina.it}   
\bigskip \\
Department of Electrical Engineering and Information Technologies, \\University of Naples Federico II, I-80125 Naples, Italy.
}

\usepackage{amsmath}
\usepackage{graphicx}
\begin{document}

\maketitle

\begin{abstract}
The aim of the paper is to present a novel class of time-dependent controls to realize ultra-fast magnetization switching in nanomagnets driven by spin-torques produced by spin-polarized electric currents. Magnetization dynamics in such systems is governed by the Landau-Lifshitz-Slonczewski equation which describes the precessional motion of (dimensionless) magnetization vector on the unit-sphere. The relevant case of nanoparticles with uniaxial anisotropy having in-plane easy and intermediate axes and out-of-plane hard axis is considered. By exploiting the characteristic smallness of damping and spin-torque intensity, the aforementioned controls are constructed via suitable perturbative tools in a way to realise approximate \emph{latitudinal solutions} (i.e.  motions on a sphere in which the out-of-plane magnetization component stays constant) with the effect to fast ``switch'' the system from one stationary state to another. The possibility to keep a (``small'') bounded value of the out-of-plane coordinate throughout this process of ``transfer'', turns out to be advantageous in the applications as it sensibly reduces the post-switching relaxation oscillations that may cause the failure of switching in real samples. Further relevant quantitative results on the behaviour of the solutions during the pre- and post- switching stages (termed ``expulsion'' and ``attraction'', respectively), are given as a byproduct. A selection of validating numerical experiments is presented alongside the corresponding theoretical results.

 
\end{abstract}

\begin{figure}[t]
    \centering
    \includegraphics[width=7cm]{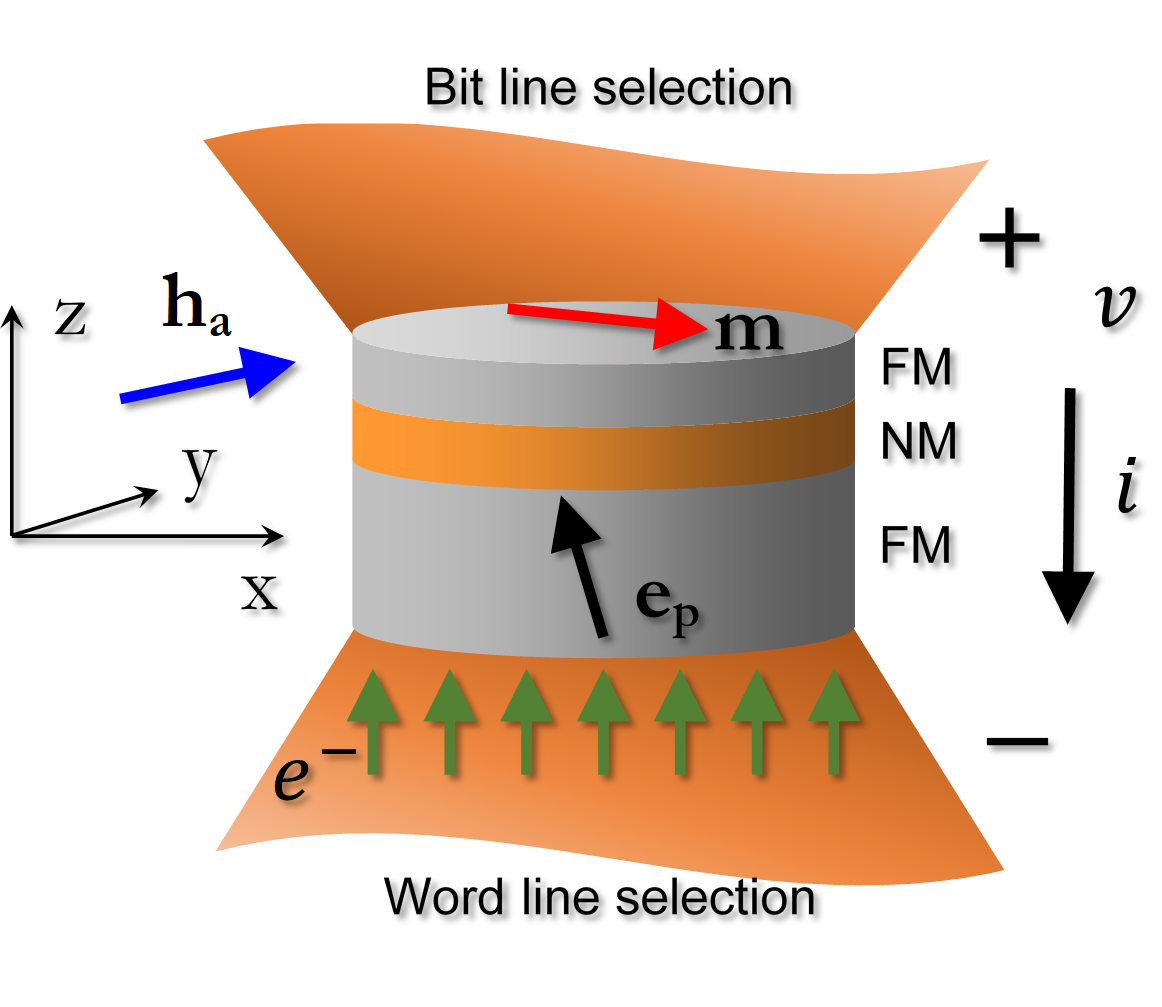}
    \caption{Sketch of a MTJ used as STT-MRAM cell.}
    \label{fig:MTJ sketch}
\end{figure}

\section{Introduction}

The efficient and high-speed manipulation of magnetic nanoelements holds immense significance in the framework of magnetization dynamics, particularly in the context of magnetic storage nanodevices and spintronics. Over the past decade, an extensive body of research has been dedicated to investigating ultra-fast magnetization switching within spintronic devices\cite{dieny2020opportunities}, which hold promise as potential candidates for advancing Magnetic Random Access Memory (MRAM) technology \cite{apalkov2016magnetoresistive,bhatti2017spintronics}.

In the pursuit of achieving rapid magnetization switching, researchers have turned their attention to the utilization of electric current pulses via the spin-transfer torque (STT) effect \cite{slonczewski1996current}. This approach has shown the potential to facilitate sub-nanosecond switching in devices like MRAM cells \cite{bedau2010spin}. Nevertheless, it is worth noting that ballistic spin-torque magnetization switching engendered by pulsed injected currents encounters a challenge reminiscent of the circumstance faced in precessional switching via transverse external magnetic fields \cite{kaka2002precessional}.

Indeed, both of these switching methodologies demand meticulous synchronization of the excitation pulse timing to ensure effective switching occurs at the precise moment\cite{daquino2004numerical}. This critical timing requirement is essential for avoiding unsuccessful switching and ensuring reliable performance. Once the excitation pulse is appropriately timed, the magnetization state transits from a high-energy to a low-energy configuration, ultimately reaching an equilibrium magnetization state \cite{bertotti2003geometrical}. Intriguingly, the underlying relaxation mechanism driving this transition is inherently stochastic, even in scenarios where thermal fluctuations are relegated to a secondary role \cite{serpico2009analytical}.

This inherent stochasticity can be attributed to the intricate interplay between the extreme sensitivity of the system to initial conditions and its multistable nature \cite{bertotti2013probabilistic}. These factors collectively bestow upon the system a probabilistic dimension, contributing to the inherent randomness in the relaxation process. This intrinsic unpredictability presents both challenges and opportunities for controlling and optimizing spin torque-induced magnetization switching in nanoscale devices.

In this respect, there is great interest in developing strategies to achieve fast, reliable and energy-efficient spin-torque switching. It has been shown that the reliability of spin-torque ballistic switching in in-plane ferromagnetic nanodots can be strongly improved by controlling the quasi-randomness by using appropriate bias fields \cite{daquino2015analysis}. On the other hand, optimization of current pulse design targeted at minimization of the energy cost has been recently proposed \cite{vlasov2022optimal}. 

In this paper, we propose a full analytical treatment for precessional spin-torque switching of in-plane magnetized nanomagnets that relies on the idea of forcing magnetization to evolve on controlled quasi-latitudinal (CQL) trajectories associated with lower values of the ferromagnet's free energy, which drive magnetization very close to the target reversed state. This is achieved by exciting magnetization precession using injected current pulses of suitable shape. This switching scheme minimizes ringing phenomena that are the main source for failure of the switching process due to back-hopping into the original magnetization state \cite{min2009back,devolder2020back}. By developing suitable perturbation theory, the rigorous conditions for the realization of CQL magnetization switching as well as the analytical expression for the optimal current pulse shape are derived as function of the system physical parameters and summarised in dedicated Lemmata.

The paper is organized as follows. In section \ref{sec:LLS dynamics}, a simple schematic of the magnetic nanosystem studied as the archetypal for in-plane STT-MRAM cells is introduced and the equation governing the magnetization dynamics driven by the spin-transfer torque is recalled along with the relevant parameters of the system. Then, after some preliminaries in sec. \ref{sec:one}, the main result of the paper is stated in sec. \ref{sec:two} as Theorem \ref{thm:main}. The development of this result requires the definition and the detailed study of three stages of the spin-torque precessional switching, termed expulsion, transfer and attraction, which is carried out in the sections \ref{sec:three}-\ref{sec:five}. Finally, validation of the proposed approach that demonstrated its effectiveness is presented.

\section{Magnetization dynamics driven by Spin-Transfer Torque}\label{sec:LLS dynamics}

We consider the magnetic nanosystem sketched in fig.\ref{fig:MTJ sketch}, representing a magnetic tunnel junction (MTJ) sandwiched between two electrodes and subject to a current of intensity $i$ flowing through the MTJ in the $z$ direction. A typical MTJ structure \cite{hirohata2020review} is composed by two ferromagnetic layers separated by a non-magnetic (NM) (insulating) layer. The magnetization of one layer, termed fixed (lower layer labelled as FM in fig.\ref{fig:MTJ sketch}), is artificially pinned to a given orientation and acts as polarizer for electron spins flowing through it. The second ferromagnetic layer, termed free (upper FM layer in fig.\ref{fig:MTJ sketch}), is where magnetization dynamics can take place driven by the external actions (injected current and applied magnetic field). The MTJ acts as a single bit-cell of a MRAM cell array where the bit state $0,1$ is coded into the mutual orientation of the free and fixed layers, namely parallel (P) or anti-parallel (AP) magnetization, respectively. The switching of the bit is triggered by the spin-transfer (STT) torque created by the electric current produced when both the bit and word lines are simultaneously addressed \cite{apalkov2016magnetoresistive,bhatti2017spintronics}, which is able to eventually switch the magnetization in the free layer of the MTJ.   

The mechanism governing the spin-torque switching can
be conveniently explained and quantified under the assumption of
spatially-uniform magnetization and negligible thermal fluctuations. In this situation, magnetization dynamics is
described by the Landau-Lifshitz-Slonczewski (LLS) equation
(written in dimensionless form)\cite{slonczewski1996current,mayergoyz2009nonlinear}:
\begin{equation}\label{eq:LLS}
\frac{d \bm m}{dt}=
-\bm m\times\bm h_\mathrm{eff}-\alpha\bm m \times(\bm m\times\bm h_\mathrm{eff})+\beta\bm m\times(\bm m\times\bm e_p)
\quad,
\end{equation}
where $\bm m$ is the magnetization unit-vector, $\bm h_\mathrm{eff}$ is the
effective magnetic field normalized by the saturation
magnetization $M_s$, $\alpha$ is the Gilbert damping, $\beta$ is
the normalized current which measures the strength of the
spin-transfer torque, $\bm e_p$ is the fixed layer unit-vector
acting as polarizer for the traversing electric current, and time
is measured in units of $(\gamma M_s)^{-1}$ ($\gamma$ is the
absolute value of the gyromagnetic ratio). The function $\beta(t)$
is proportional to the injected current $i(t)$ through the
relationship $\beta(t)=b_p i(t)/(S J_p)$  where $b_p$ is a
model-dependent parameter in the order of unity,
 $S$ is the device cross sectional area, and $J_p=\mu_0
M_{s}^{2} |e|d/\hbar$ is a characteristic current density ($\mu_0$
is the vacuum permeability, $e$ is the electron charge, $d$ is the
thickness of the free layer, and $\hbar$ is the reduced Planck
constant).

The effective field is expressed as:
\begin{equation}\label{eq:effective field}
\bm h_\mathrm{eff}=-\frac{\partial g_L}{\partial\bm m}=
-D_1 m_1\bm e_1-D_2 m_2\bm e_2-D_3 m_3\bm e_3+ \bm h_a \quad,
\end{equation}
where $\bm e_1,\bm e_2,\bm e_3$ are the cartesian unit-vectors
along the coordinate axes, $g_L$ is the magnetic free
energy of the particle
\begin{equation}\label{eq:free energy}
g_L(\bm m,\bm h_a)=\frac{1}{2}D_1 m_1^2+\frac{1}{2}D_2
m_2^2+\frac{1}{2}D_3 m_3^2-\bm h_a\cdot\bm m \,\,,
\end{equation}
$\bm h_a$ is the external applied magnetic field and $D_1,D_2,D_3$
are effective demagnetizing factors taking into account shape and
magneto-crystalline anisotropy. In the sequel, we assume
\beq{eq:di}
0<D_1<D_2<D_3 \mx{,}
\eeq
meaning that $\bm e_1,\bm e_2,\bm e_3$ refer to the easy,
intermediate, hard axes, respectively.

In the absence of injected current $\beta=0$ and small enough bias applied field $\bm h_a$, magnetization lies in one of the stable equilibria aligned with the easy axis. We assume that this equilibrium is $s^-$ such that $m_1\approx -1$.  
Under this assumption, it is apparent from eq.\eqref{eq:LLS} that, in order to have the maximum torque on magnetization when the current is turned on, the polarizer $\bm e_p$ must be orthogonal to the easy axis. In this paper we assume $\bm e_p=\bm e_3$ for reasons that will be clear in the sequel.

\section{Preliminaries and set-up}\label{sec:one}
Leu us consider the following system of ODEs 
\beq{eq:prsw}
\begin{cases}
\dot{u}_1=h_2 u_3 + D_{3,2} u_2 u_3 +  u_1 u_3 \beta(t) + \alpha u_1 \left(D_{3,1} u_3^2 + D_{2,1} u_2^2 - h_2 u_2 \right) \\
\dot{u}_2=-D_{3,1} u_1 u_3 + u_2 u_3 \beta(t) + \alpha \left( D_{3,2} u_2 u_3^2 - D_{2,1} u_2 u_1^2 +h_2 u_1^2+h_2 u_3^2\right) \\
\dot{u}_3=-h_2 u_1 + D_{2,1} u_1 u_2 -\left(u_1^2+u_2^2 \right) \beta(t) -\alpha u_3 \left( D_{3,2} y^2 +D_{3,1} u_1^2 +h_2 u_2 \right)
\end{cases}
\mx{.}
\eeq
The latter is immediately obtained from (\ref{eq:LLS}) by setting $\ue{u}:=\ue{m}$ and $\ue{h}_a=:(0,h_2,0)$. Given $D_i$ satisfying (\ref{eq:di}), we have defined, for all $ 1 \leq j < i \leq 3$,  
\[
D_{i,j}:=D_i-D_j \mx{.}
\]
Clearly, these quantities are such that
\beq{eq:propertiesd}
D_{3,1}>D_{3,2}, \qquad D_{3,2}+D_{2,1}-D_{3,1}=0 \mx{.}
\eeq
The function $\beta(t):[0,+\infty) \rw \RR$ (\emph{injected current}) plays here the role of ``control''.\\
As it is well known (see, e.g. \cite{mayergoyz2009nonlinear}), the function ``distance from the origin'' 
\beq{eq:energy}
\Psi(\ue{u}):=\sum_{j=1}^3 u_j^2
\eeq
is conserved along the solutions of (\ref{eq:prsw}), namely $\Psi(\ue{u})=1$, meaning that the motion occurs on the unit-sphere. \\
A switching process (see \cite{mayergoyz2009nonlinear} for a comprehensive treatment) can be briefly described as the  possibility to ``move'' the state of the system (\ref{eq:prsw}) from (a neighbourhood of) the equilibrium $\ue{s}^{-}$ to  $\ue{s}^{+}$, where $\ue{s}^{\pm}:=(\pm \sqrt{1-(h_2/D_{2,1})^2}, h_2/D_{2,1},0)$ when $\beta(t)=0$ (straightforward check) and this can be achieved by applying a suitable injected current $\beta(t)$. We shall say that the ``switching'' has taken place, once the system has reached a state which is attracted by $\ue{s}^{+}$ without the need to supply any further control via $\beta(t)$. For reasons of topological nature, such a transition cannot exist if the current is switched off all the time.  \\
As already mentioned in the foreword, several approaches can be found in the literature. Possibly the most paradigmatic one is the ballistic one, which  consists in applying a current $\beta(t)=\const >0$ for a prescribed period of time $t \in [0,T_e]$, then switching it off. The time $T_e$ is meant to be be determined. \\ 
However, as already stressed, due to the highly involved phase space space structure of the system at hand pointed out by specialised analyses (see, e.g., \cite{10.1063/1.4914122}) such as  a successful switching is more properly intended as a probabilistic feature of the described technique. In other terms, without any supporting result of quantitative nature, this class of switching approaches have an experimental validity only. \\
In a different spirit, this paper aims to provide a \emph{deterministic} argument to successfully complete a switching process under mild assumptions. This has been possible by recognising some intrinsic perturbative features of the system at hand, and exploiting them via perturbation methods borrowed from the Hamiltonian world. As a main achievement, this has led to the realisation of the mentioned CQL solutions, via a suitably constructed controlling current $\beta(t)$. As the name suggests, these solutions are characterised by a quasi-constant (in the perturbative sense) $u_3(t)$ throughout the motion, until a neighbourhood of the target equilibrium $\ue{s}^{+}$ is reached.\\  Interestingly enough, ``small'' values of $|u_3(t)|$ turn out to provide a very adequate option, being the solution ``rapidly'' attracted by the equilibrium (characterised by $u_3=0$), right after the current is switched off. This constitutes another difference with the ballistic approach, where the required amount of initial injected current leads to ``higher'' values of $|u_3(t)|$, requiring in this way a ``long time'' to come back to zero, as per effect of the ``small'' dissipation $\alpha$. 

\section{Perturbative setting and main result}\label{sec:two}
The starting point consists in noticing that a class of realistic models exhibit different scales amongst the involved parameters. More specifically, whilst $D_{3,2},D_{3,1}$ are ``of order one'', denoted with $O(1)$, the remaining $h_2,D_{2,1},\alpha$ are sensibly smaller, namely ``of order $\lambda$'', where $\lambda$ is typically $\sim 10^{-2}$. We choose $\beta(t) = O(\lambda)$ as it is usually done, since the spin torque compensates the damping and is of the same order of magnitude. The described feature of the parameters at hand, naturally leads to a perturbative formulation for (\ref{eq:prsw}). In particular, by defining $h_{2}=:\lambda \tilde{h}_2$, with $\tilde{h}_2=O(1)$ and similarly for all the other mentioned $O(\lambda)$ quantities, the system (\ref{eq:prsw}) is immediately cast into the following form
\beq{eq:prswscaled}
\begin{cases}
\dot{u}_1=D_{3,2} u_2 u_3 +  \lambda \left[ \thd u_3+ u_1 u_3 \tb (t) + \ta u_1 D_{3,1} u_3^2 \right]+\ta \lambda^2 u_1 u_2 \left[  \tilde{D}_{2,1} u_2 - \thd   \right]\\
\dot{u}_2=-D_{3,1} u_1 u_3 +  \lambda u_2 u_3 \left[  \tb (t) + \ta D_{3,2} u_3 \right] 
+\ta \lambda^2 \left[\thd \left( u_1^2+ u_3^2 \right) -\tilde{D}_{2,1}  u_1^2 u_2 \right] \\
\dot{u}_3=\lambda \left[-\thd u_1 + u_1 u_2 \tilde{D}_{2,1} -\left(u_1^2+u_2^2 \right) \tb (t)  
-\ta u_3 \left( D_{3,2} u_2^2 +D_{3,1} u_1^2 \right) \right] -\ta \thd \lambda^2 u_2 u_3
\end{cases}
\eeq
or, in a more compact notation,
\beq{eq:prswscaledalternative}
\dot{\ue{u}}=\uue{A}(u_3)\ue{u}+\lambda \ue{\mathcal{N}}(\ue{u},\lambda;\tb (t)) \mx{,}  
\eeq 
where
\[
\uue{A}(\xi):=\begin{pmatrix}
	0 & D_{3,2}\xi & 0\\
	-D_{3,1}\xi & 0 & 0\\
	0 & 0 & 0\\
\end{pmatrix}
\]
and the definition of $\ue{\mathcal{N}}$ is obvious.\\
It is clear that at the order zero in $\lambda$, the variable $u_3(t)$ is constant for all $t$. As a consequence, the dynamics in the variables $(u_1,u_2)$ is described by the equation of a harmonic oscillator with (constant) frequency determined once and for all by the initial condition $u_3(0)$. We shall restrict ourselves to the choice $u_3(0)=-K$, with $K \in \RR^+$. We remark that the latter will be achieved in the full process via the ``expulsion'' stage, which will be described later on.\\
Hence, by defining
\beq{eq:omega}
\omega:=K \sqrt{D_{3,2} D_{3,1}}
\eeq
the energy of such an oscillator is given by $E=(\omega^2/2)[u_1^2+(u_2/D_{3,1})^2]$, and hence a family of ellipses in the plane $(u_1,u_2)$ parameterised by $E$. A motion starting in $(u_1(0),u_2(0))=(-1,0)$ (i.e. $\lambda$-close to $(s_1^-,s_2^-)$) has energy $\hat{E}=\omega^2/2$. As a consequence, after a semi-period 
\beq{eq:ttrplus}
T_{tr}^+:=\pi/\omega \mx{,}
\eeq
it will evolve in $(u_1(T_{tr}^+),u_2(T_{tr}^+))=(1,0)$, which is $\lambda$-close to $(s_1^+,s_2^+)$. This is, essentially, the key mechanism of transfer we will rely on.\\
Note that this represents an archetypal \emph{latitudinal} solution, as (at least at the zero-th order in $\lambda$) it takes place along the same latitudinal line of the unit sphere, determined once and for all by the choice of $K$. \\
However, should $K$ be ``small'' (for reasons that will be clarified later on), the transfer process would take a ``long'' time, and then the $O(\lambda)$ contributions would become significant. This means that, in order to preserve the ``transfer'' feature of the solutions, we need to find a way to ``contain'' the variations of $\ue{u}$ from the unperturbed solution. \\
At this point, we are about to summarise the argument used to construct the control that is apt to realise the desired switching process.\\
Let us now consider, for this purpose, the auxiliary system in a renamed set of variables, say $v_j(t)$, which is obtained from (\ref{eq:prswscaled}) when $\tb(t)$ is replaced by some function of $\ue{v}$, namely $\hat{\beta}(\ue{v})$ and subject to the same initial conditions. Hence, the key step to obtain a latitudinal solution consists in choosing $\hat{\beta}$ in such a way the r.h.s. of the third equation is zero. In fact, as a consequence of this choice, we have $v_3(t)=-K$ for all $t$. \\
Hence, in order to select a class of latitudinal solutions we can choose  
\beq{eq:choicebeta}
\hat{\beta}_{lat}(\ue{v}):=(v_1^2+v_2^2)^{-1} \left[ v_1 v_2 \tilde{D}_{2,1} + \ta K \left( D_{3,2} v_2^2 +D_{3,1} v_1^2 \right) + \ta \thd \lambda v_2 K
 -\thd v_1 \right] \mx{,}
\eeq
this leads to
\beq{eq:prswscaledtwo}
\dot{\ue{v}}=\uue{A}(v_3)\ue{v}+\lambda \ue{\mathcal{N}}(\ue{v},\lambda;\hat{\beta}_{lat}(\ue{v}))
=:\uue{A}(-K)\ue{v}+\lambda \ue{\mathcal{F}}(\ue{v};\ue{v})+\lambda^2 \ue{\mathcal{R}}(\ue{v};\ue{v}) \mx{.}
\eeq
In other terms,  $\lambda \ue{\mathcal{F}}(\ue{v};\ue{v})$ denotes the first order terms (in $\lambda$) of the field $\ue{\mathcal{N}}(\ue{v},\cdot;\cdot)$ in which $\hat{\beta}$ has been substituted as a function of $\ue{v}$. It might be useful, in order to clarify the notation, to denote as $\lambda \ue{\mathcal{F}}(\ue{u};\ue{v})$ the first order terms arising in (\ref{eq:prswscaledalternative}) following to the substitution $\tilde{\beta} \leftarrow \hat{\beta}_{lat}(\ue{v})$. The meaning of $\ue{\mathcal{R}}$ is obvious as a consequence. Their explicit expression is slightly more involved and is given in Appendix A, however, as for $\ue{\mathcal{F}}(\ue{v};\ue{v})$ and $\ue{\mathcal{R}}(\ue{v};\ue{v})$, we have  
\beq{eq:remaindersone}
\begin{aligned}
\mathcal{F}_1(\ue{v};\ue{v}) & =-\rho K v_2 \left(\thd v_2 + \tilde{D}_{2,1} v_1^2 \right) \\
\mathcal{F}_2(\ue{v};\ue{v}) & = \rho K v_1 \left( \thd v_2 - \tilde{D}_{2,1} v_2^2 \right) \\
\mathcal{R}_1(\ue{v};\ue{v}) & = \rho \ta v_1 v_2 \left[ K^2 \tilde{D}_{2,1} v_2 + v_1^2 (v_2 \tilde{D}_{2,1}-\thd) -  \thd (K^2 + v_2^2) + \tilde{D}_{2,1} v_2^3 \right]\\
\mathcal{R}_2(\ue{v};\ue{v}) & = \rho \ta v_1^2 \left[ - K^2 \tilde{D}_{2,1} v_2 - v_1^3 (v_2 \tilde{D}_{2,1}-\thd) +  \thd (K^2 + v_2^2) - \tilde{D}_{2,1} v_2^3 \right]
\end{aligned} \mx{,}
\eeq
where, by virtue of (\ref{eq:energy}), we have set $v_1^2+v_2^2=1-K^2=:1/\rho$. Clearly, by construction 
\beq{eq:remaindersonethree}
\mathcal{F}_3(\ue{v};\ue{v})=\mathcal{R}_3(\ue{v};\ue{v})=0 \mx{.}
\eeq
Hence (\ref{eq:prswscaledtwo}) is, by all means, a two dimensional system of ODEs.\\
Similarly to what has been observed for (\ref{eq:prswscaled}), (\ref{eq:prswscaledtwo}) is clearly integrable at the zero-th order in $\lambda$. \\
Let us imagine for a moment to be able to compute a solution $v_{1,2}(t)$ to (\ref{eq:prswscaledtwo}). The construction above implies that, by choosing the particular \emph{control} 
\beq{eq:betachoice}
\tilde{\beta}(t):= \hat{\beta}_{lat}(v_{1}(t),v_{2}(t),-K) \mx{,}
\eeq
and substituting it in (\ref{eq:prswscaled}), one has $\ue{u}(t)=(v_{1}(t),v_{2}(t),-K)$, provided that (\ref{eq:prswscaled}) and  (\ref{eq:prswscaledtwo}) are subject to the same initial conditions. In other terms, (exact) solutions to (\ref{eq:prswscaledtwo}) yield latitudinal solutions (in the sense defined above) to (\ref{eq:prswscaled}) via (\ref{eq:betachoice}). Hence, it is reasonable to expect that ``approximated'' solutions to (\ref{eq:prswscaledtwo}) could provide a ``quasi-constant'' behaviour for the third variable $u_3(t)$, in a sense to be made precise. We shall refer to these as \emph{controlled quasi-latitudinal} (CQL) solutions to (\ref{eq:prswscaled}). Under this setting, the main result states as follows:
\begin{satz}\label{thm:main} 
There exist $\lambda_0,r^-,T_e,T_{tr}>0$, a non-trivial set of parameters and a control 
\beq{eq:control}
\tilde{\beta}(t):=
\begin{cases}
\tilde{\beta}_e & t\in [-T_e,0)\\
\tilde{\beta}_{tr} (t) & t \in [0,T_{tr}]\\
0 & t \in (T_{tr}, +\infty)
\end{cases}
\eeq
such that the flow of the system (\ref{eq:prswscaled}) controlled via (\ref{eq:control}) and denoted with $\tilde{\Phi}^t$, for any $r\leq r^+$ and initial condition $\ue{u}_0:=\ue{u}(-T_e) \in \mathfrak{B}_{r} ({\ue{s}^-})$, satisfies 
\beq{eq:transferstat}
\lim_{t \rw +\infty} \tilde{\Phi}^t (\ue{u}_0) \in \mathfrak{B}_{f r} ({\ue{s}^+})  \mx{,}
\eeq
for some $f>0$. In particular, the system realises a CQL solution in the time interval $[0,T_{tr}]$.
\end{satz}
The rest of the paper is devoted to the proof of Thm. \ref{thm:main}. This will be achieved in three steps. The first one, called ``expulsion'', consists in choosing a suitable neighbourhood of the point $\ue{s}^-$ and injecting a constant current of magnitude $\tilde{\beta}_e$ for a certain time $T_e$. It is shown how this ``pushes'' the third variable $u_3(t)$ in a region in which it is strictly negative, a key requirement in order to obtain the class of switching motions we are interested in, even at the zero-th order, as anticipated in sec. \ref{sec:two}.\\
A second stage, referred to as ``transfer'', represents the very heart of the argument, i.e. where the QLS are constructed via a suitable controlling current, with the property of being able to ``deliver'' the solutions starting in the vicinity of the point reached once the first step is complete, in a neighbourhood of the target point $\ue{s}^+$. \\
The last step, called ``attraction'', consists in showing that $\ue{s}^+$ is an attractive point when the current is switched off, attracting in this way any points arriving in its vicinity. Although this is a very well known property, quantitative information are needed for a precise formulation of the result, and this is exactly the aim of this third stage. It is immediate to realise that being able to carry out the described three steps, implies a successful realisation of the switching process. The whole argument clearly exploits the known group property of a dynamical system. \\
From a technical viewpoint, the ``core'' of the proof (consisting in the ``transfer'' stage), relies on the possibility to construct a suitable control for the system at hand, by approximating the non-linear flow via Hamiltonian perturbative tools. This is possible via the well known possibility to interpret any system of ODEs as a Hamiltonian system in a suitably extended phase space, see e.g. \cite{berdvar}. However, such a simple observation has been profitably used in several cases, due to the potential of a full all-orders generalisation of the perturbative setting. See, for instance, \cite{for162}, \cite{fliap}. This has the potential to increase the threshold $\lambda_0$, expanding in this way the class of systems which can be dealt with via this approach.\\
The remaining stages, ``expulsion'' and ``attraction'', are carried out by using established tools from the theory of ODEs that will be specified later.  


\section{Expulsion}\label{sec:three}
\begin{figure*}[t!]\begin{center}
		\begin{minipage}[c][1\width]{0.48\textwidth}
			\hspace{0pt} 
			{\begin{overpic}[width=\textwidth]{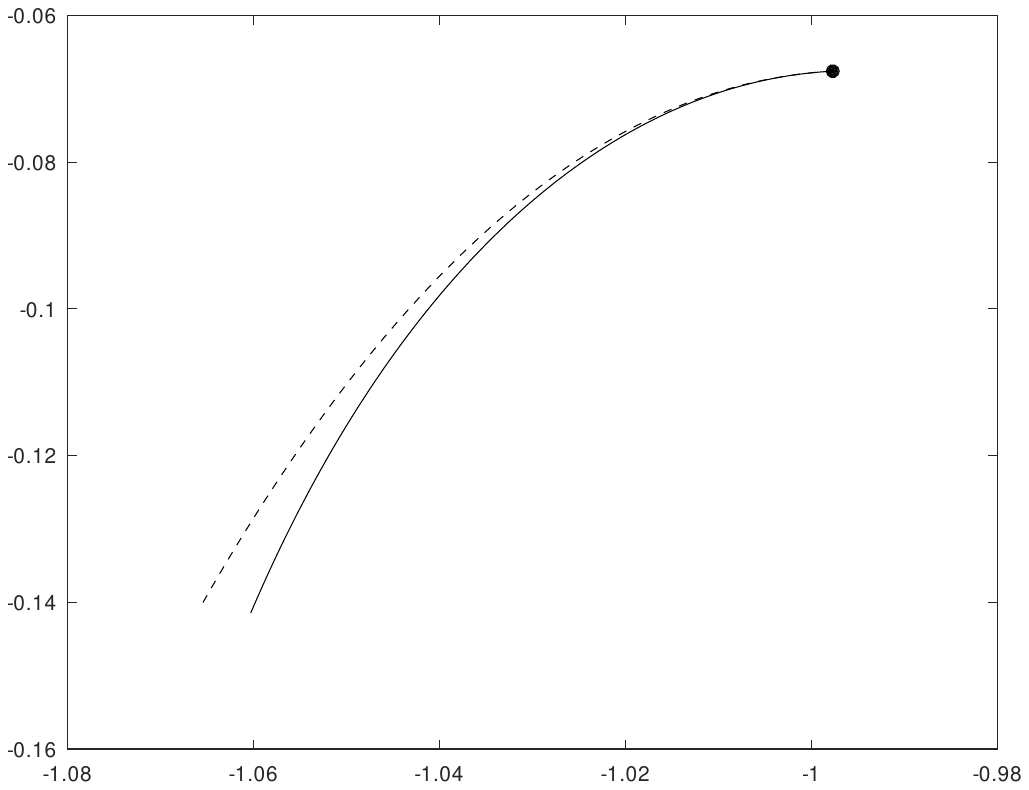}
					\put(0,78){\footnotesize (a)}
	                \put(5,68){\footnotesize $w_1$}
					\put(60,24){\footnotesize $w_2$}
                \end{overpic}}
		\end{minipage}
    		\begin{minipage}[c][1\width]{0.48\textwidth}
			\hspace{-5pt}
			{\begin{overpic}[width=1\textwidth]{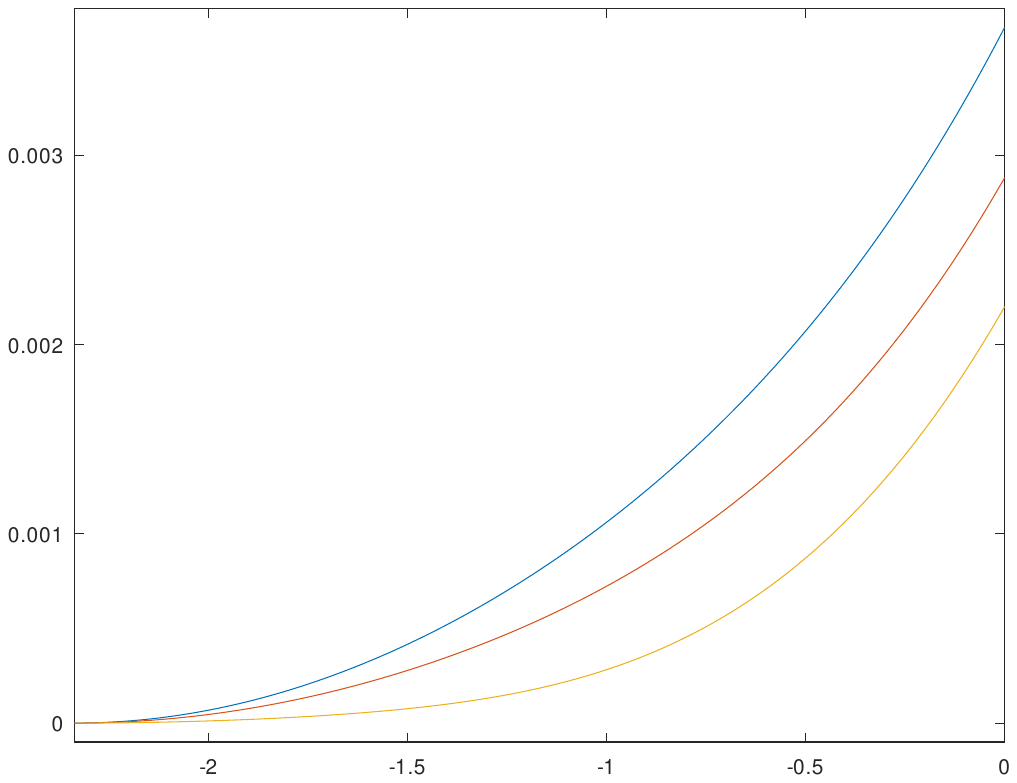}
					\put(0,78){\footnotesize (b)}
			 	  \put(5,68){\footnotesize $|\ue{\delta}|$}
					\put(62,24){\footnotesize $t$}
					\put(70,71){\footnotesize $\lambda_{(1)}$}
   					\put(70,63){\footnotesize $\lambda_{(2)}$}
   					\put(70,55){\footnotesize $\lambda_{(3)}$}
			\end{overpic}}
		\end{minipage}	\\
		\vspace{-30pt}
  \end{center}
	\caption{In panel (a) the trajectories $\ue{\xi}(\tau)$ (continuous line) and $\ue{\xi}'(\tau)$ (dashed line) are reported in the stereographic coordinates defined via $w_{1,2}:=u_{1,2}/(1+u_3)$ and transformation (\ref{eq:translation}). Panel (b) shows the behaviour of the error $|\ue{\delta}(t)|$ for $t \in [-T_e,0]$ with $T_e=2.3372$ and three different values of $\lambda$: $\lambda_{(1)}=0.006$, $\lambda_{(2)}=0.004$ and $\lambda_{(3)}=0.001$. The corresponding values for $\tilde{\beta}_e$ are $0.0018$, $0.0012$ and $0.003$, respectively. The remaining parameters used are specified in Appendix C.}
	\label{fig:expulsion}
	\vspace{0pt}
\end{figure*}

Let us consider the system (\ref{eq:prswscaledalternative}) and set $\thd=:-\tilde{D}_{2,1}\Omega$, with $\Omega>0$. In this setting, in particular, one has $\ue{s}^{\pm} \equiv (\pm \gamma,-\Omega,0)$, where $\gamma:=\sqrt{1-\Omega^2}$. Let us now perform a translation-rescaling of $\hat{\ue{s}}^{-}$ at the origin of the new system of coordinates 
\beq{eq:translation}
\lambda \ue{\xi}:=\ue{u}-\ue{s}^{-} \mx{.}
\eeq  
Throughout this section we shall consider the autonomous control $\lambda \tilde{\beta}(t) =: \lambda \tilde{\beta}_e(t) \equiv {\beta}_e=\const$ and the new time $\tau:=T_e+t$, for any $T_e>0$. Consequently, $\Phi_e^{\tau}$ will denote the phase flow of (\ref{eq:prswscaledalternative}) with respect to the new time. Hence, by (\ref{eq:translation}), (\ref{eq:prswscaledalternative}) reads as 
\beq{eq:systemexp}
\frac{d}{d \tau} \ue{\xi}=
\uue{\ml{L}}\ue{\xi}+\ue{f}+\lambda \ue{\ml{V}}(\ue{\xi}; \lambda,\beta_e) \mx{,}
\eeq
where, 
\[
\uue{\ml{L}}:=\begin{pmatrix}
	0 & 0 & \bar{a} \\
	0 & 0 & \bar{b} \\
	2 \beta_e \gamma & 2 \beta_e \Omega & 0\\
\end{pmatrix} 
\mx{,}
\qquad
\ue{f}=(0,0,-\beta_e /\lambda)^{\top} \mx{,}
\]
\beq{eq:abarbbar}
\bar{a}:=-D_{3,2} \Omega-\beta_e \gamma, \qquad \bar{b}:=D_{3,2} \gamma-\beta_e \Omega \mx{,}
\eeq and $\ue{\ml{V}}(\ue{\xi}; \lambda,\beta_e)$ is defined as a consequence. Let us remark that we have used the relation $D_{3,1}=D_{3,2}-\lambda \tilde{D}_{2,1}$, see (\ref{eq:propertiesd}), in order to simplify the structure of $\uue{\ml{L}}$. As a result, the linear term $-\tilde{D}_{2,1}\xi_3$ is included in  $\ue{\ml{V}}(\ue{\xi}; \lambda,\beta_e)$. In particular, it is evident that for sufficiently small $\lambda$ the dominant behaviour is ``expulsive'' i.e. it drives $\xi_3$ towards negative values.  \\
As it can be easily checked, by defining 
\[
\uue{S}:=\begin{pmatrix}
	1 & 1 & 1\\
	\bar{b}/\bar{a} & \bar{b}/\bar{a} & -\gamma/\Omega\\
	\sqrt{2} i \beta_e /\bar{a} & -\sqrt{2} i \beta_e / \bar{a} & 0\\
\end{pmatrix} 
\]
one has $ \uue{S}^{-1} \uue{\ml{L}} \uue{S}= \uue{\Gamma} \equiv \diag (-\sqrt{2} i \beta_e, \sqrt{2} i \beta_e,0)$. Hence, the solution to (\ref{eq:prswscaledalternative}) reads as 
\beq{eq:solutiontoexp}
\ue{\xi}(\tau)=e^{\uue{\ml{L}} \tau} \ue{\xi}(0)+\int_0^{T_e} e^{\uue{\ml{L}} (\tau-s)} \ue{f} ds + 
\lambda \int_0^{T_e} e^{\uue{\ml{L}} (\tau-s)} \ue{\ml{V}}(\ue{\xi}; \lambda,\beta_e) ds  \mx{,}
\eeq
where, by noticing that $\bar{a}\gamma+\bar{b}\Omega=-\beta_e$, one has
\[
e^{\uue{\ml{L}} \tau}=\frac{1}{\beta_e}
\begin{pmatrix}
-\bar{b}\Omega - \bar{a} \gamma \cos (\sqrt{2} \tau \beta_e)& 
 \bar{a}\Omega - \bar{a} \Omega \cos (\sqrt{2} \tau \beta_e)& 
(\bar{a}/\sqrt{2}) \sin (\sqrt{2} \tau \beta_e)
\\
\bar{b} \gamma - \bar{b} \gamma \cos (\sqrt{2} \tau \beta_e)   & -\bar{b} \Omega \cos (\sqrt{2} \tau \beta_e) -\bar{a} \gamma  &  (\bar{b}/\sqrt{2}) \sin (\sqrt{2} \tau \beta_e) \\
\sqrt{2} \beta_e \gamma \sin (\sqrt{2} \tau \beta_e) & \sqrt{2} \beta_e \Omega \sin (\sqrt{2} \tau \beta_e)  & \beta_e \sin (\sqrt{2} \tau \beta_e) \\
\end{pmatrix} \mx{.}
\]
Furthermore, the norm of the latter is uniformly bounded in time, i.e.
\beq{eq:expbound}
|| e^{\uue{\ml{L}} \tau} || \leq \ml{M}_e \mx{,}
\eeq
for some suitable $\ml{M}_{e}=\ml{M}_e(\beta_e,D_{3,2},\gamma,\Omega) \geq 1$ and all $\tau \in \RR$. \\
It is now easy to realise that, if $O(\lambda)$ are disregarded in (\ref{eq:solutiontoexp}), the integral appearing in it can be computed immediately and this gives rise to a function, namely $\ue{\xi}'(\tau)$, which approximates the ``full'' $\ue{\xi}(\tau)$ for ``small'' $\lambda$ (under suitable assumptions). Such a function reads as 
\beq{eq:apprxi}
\ue{\xi}'(\tau)=e^{\uue{\ml{L}} \tau} \ue{\xi}'(0) +(2 \beta_e \lambda)^{-1}  \left(\bar{a}(\cos (\sqrt{2}t \beta_e)-1),\bar{b}(\cos (\sqrt{2}t \beta_e)-1),-\sqrt{2} \beta_e \sin (\sqrt{2}t \beta_e) \right)^{\top} \mx{.}
\eeq
As the state of the system is supposed to be initially in a neighbourhood of the origin, it is natural to set $\ue{\xi}'(0)=\ue{0}$. Hence, in the original set of variables, the third variable evolves as $u_3(\tau)\sim -2^{-1/2} \sin (\sqrt{2}t \beta_e)$ for initial data ``close to'' $\ue{s}^-$ and ``sufficiently small'' $\lambda$. Once a target value $-K$ for $u_3$ has been set, with 
\beq{eq:existencet}
0< \sqrt{2} K \leq 1 \mx{,}
\eeq
the expulsion stage is defined as the evolution of the system for $\tau \in [0, T_e]$, where
\beq{eq:texp}
T_e:=(\sqrt{2}\beta_e)^{-1} \arcsin (\sqrt{2} K) \mx{.}
\eeq    
From (\ref{eq:apprxi}) one immediately gets
\beq{eq:apprxitwo}
\ue{\xi}'(T_e)=(2 \beta_e \lambda)^{-1} (\bar{a}(-1+\sqrt{1-2 K^2}),\bar{b}(-1+\sqrt{1-2 K^2}),-2 K \beta_e) \mx{.}
\eeq
\begin{rem}%
For realistic values of the parameters (for instance, $\beta_e,K^2=O(\lambda)$), the quantities $\lambda \xi_{1,2}'(T_e)$ are $O(1)$, i.e. the point reached after the expulsion stage is ``far'' from $\ue{s}^-$. For this reason, approximation formula (\ref{eq:apprxitwo}) will play a key role in the main proof. Furthermore, the latter has some interest itself in ballistic switching processes, in which ``long'' expulsions are typically considered.
\end{rem}
\noindent The approximation features of (\ref{eq:apprxitwo}), numerically validated for an example in fig. \ref{fig:expulsion}, are stated in a quantitative form in the next
\begin{lem}\label{lem:expulsion}
Let us set a target value $K$ satisfying (\ref{eq:existencet}) and $\beta_e>0$ arbitrarily chosen. Then choose $T_e$ as in (\ref{eq:texp}). Let us now define 
\beq{eq:rstar}
r^*:=(2 \beta_e \lambda)^{-1} \left[4 K^2 + (D_{3,2}^2+\beta_e^2)(\sqrt{1- 2 K^2}-1)^2 \right]^{\frac{1}{2}}
\eeq 
and let be $\rho_e>0$ arbitrarily chosen. Now set 
\beq{eq:lip}
\lambda \ml{M}_1:=\max_{\ue{\xi} \in \mathfrak{B}_{r^*}(\ue{0})} | \ue{\ml{V}}(\ue{\xi};\lambda,\beta_e)| \mx{,}
\qquad   
\lambda \ml{M}_2:=\max_{\ue{\xi} \in \mathfrak{B}_{r^*+\rho_e}(\ue{0})} ||\uue{D} \ue{\ml{V}}(\ue{\xi};\lambda,\beta_e)||_{\infty} \mx{,} 
\eeq
see Appendix B. for the explicit expressions of $\uue{D} \ue{\ml{V}}(\ue{\xi};\lambda,\beta_e)$.\\
Then for all $\lambda \in (0,\lambda_e]$, with
\beq{eq:lamz}
\lambda_e:=(T_e \ml{M}_e)^{-1} \min\{  (4 \mathcal{M}_1)^{-1}\rho_e, \mathcal{M}_2^{-1} \log 2 \} \mx{,}
\eeq
the following property holds
\beq{eq:shadowone}
\Phi^{T_e} (\mathfrak{B}_{(4 \ml{M}_e)^{-1} \rho_e \lambda }(\ue{s}^-)) \subseteq \mathfrak{B}_{\rho_e \lambda  }(\ue{u}_c(T_e)) \mx{,}
\eeq
where $\ue{u}_c(\tau):=\lambda \ue{\xi}'(\tau)+\ue{s}^{-}$. In other terms, solutions starting within a sphere of radius $r^{-}$ around $\ue{s}^{-}$ remain contained in a sphere of radius $ \rho_e \lambda$ of the (known) point $\lambda \ue{\xi}'(T_e)+\ue{s}^{-}$. 
\end{lem}
\begin{rem}
As it can be easily noticed from (\ref{eq:rstar}), $r^*|_{K=0}=0$ and it grows monotonically with $K$. The values of $\mathcal{M}_{1,2}$ will increase accordingly. This implies that the allowed threshold for $\lambda_e$ worsens as the target $K$ gets bigger.   
\end{rem}
\proof
The proof relies on the fact that the solutions $\ue{\xi}(\tau)$, starting from a suitable neighbourhood of the origin, do not escape the sphere of radius $r^*+\rho_e$ centred at the origin, in such a way the bounds (\ref{eq:lip}) are justified. For this purpose, let us firstly check that $|\ue{\xi}'(T_e)| \leq r^*$. This is immediate from (\ref{eq:apprxitwo}), assumption (\ref{eq:rstar}) (notice that $\bar{a}^2+\bar{b}^2=D_{3,2}^2+\Omega^2$) and finally observing that $|\ue{\xi}'(\tau)|$ is monotonically increasing for all $\tau \in [0, T_e]$. \\ 
Let us now define $\ue{\delta}(\tau):=\ue{\xi}(\tau)-\ue{\xi}'(\tau)$. By substituting in (\ref{eq:solutiontoexp}) we get
\beq{eq:soldelta}
\ue{\delta}(\tau)=e^{\uue{\ml{L}} \tau} \ue{\delta}(0)+ 
\lambda \int_0^{T_e} e^{\uue{\ml{L}} (\tau-s)} \ue{\ml{V}}(\ue{\xi}'+\ue{\delta}; \lambda,\beta_e) ds  \mx{,}
\eeq
where $\ue{\delta}(0) \equiv \ue{\xi}(0)$. By taking the absolute values of (\ref{eq:soldelta}), then using the bound $|\ue{\ml{V}}(\ue{\xi}'+\ue{\delta}; \lambda, \beta_e )| \leq |\ue{\ml{V}}(\ue{\xi}'; \lambda ,\beta_e )|+\ml{M}_2 |\ue{\delta}| $, assumptions (\ref{eq:lip}), and finally the classical Gr\"{o}nwall lemma, one obtains 
\[
|\ue{\delta}(\tau)| \leq \ml{M}_e ( |\ue{\delta}(0)| + \lambda \ml{M}_1 \tau) \exp(\lambda \ml{M}_e \ml{M}_2 \tau) \mx{.}
\] 
Hence, for all $\lambda \in (0,\lambda_e]$ and all $\tau \in (0,T_e]$ it is sufficient to choose $|\ue{\delta}(0)| \leq (4 \ml{M}_e)^{-1} \rho_e$ as suggested by the l.h.s. of (\ref{eq:shadowone}), in order to get $|\ue{\delta}(\tau)| \leq \rho_e$. This proves the r.h.s. of (\ref{eq:shadowone}). 
\endproof

\section{Transfer}\label{sec:four}
\begin{figure*}[t!]\begin{center}
  		\begin{minipage}[c][1\width]{0.48\textwidth}
			\hspace{0pt}
			{\begin{overpic}[width=1\textwidth]{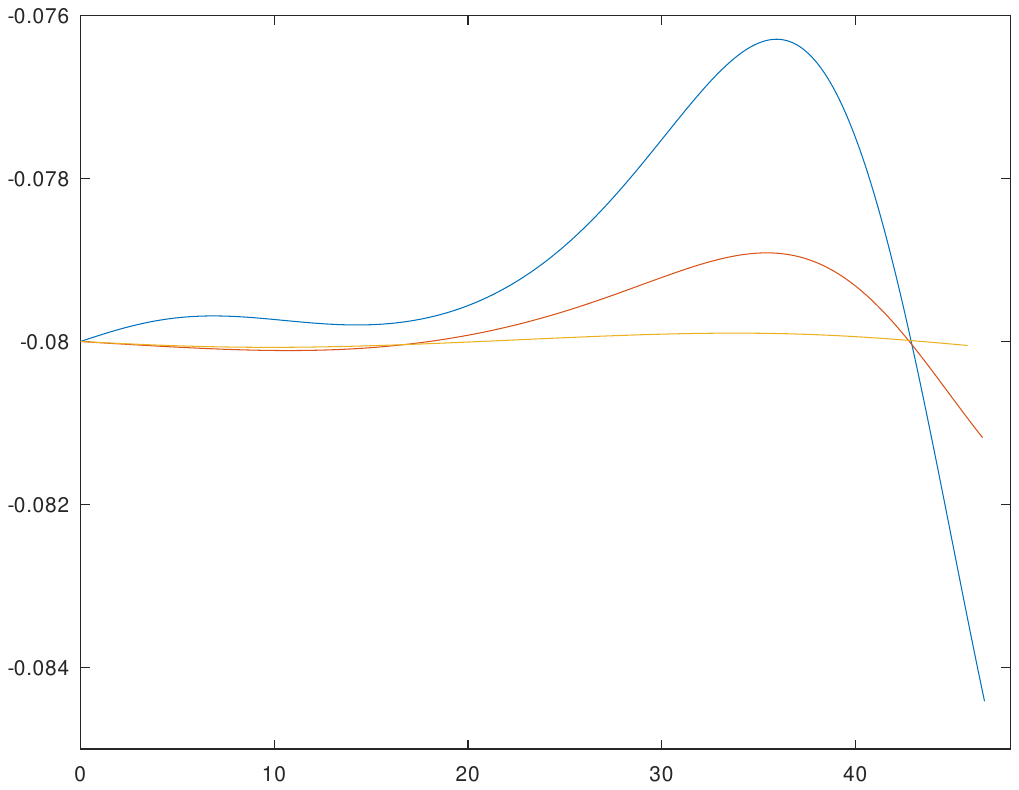}
					\put(0,78){\footnotesize (a)}
					\put(4,67.5){\footnotesize $u_3$}
					\put(63,24){\footnotesize $t$}
					\put(41,67){\footnotesize $\lambda_{(1)}$}
   					\put(50,61){\footnotesize $\lambda_{(2)}$}
   					\put(50,50.5){\footnotesize $\lambda_{(3)}$}
			\end{overpic}}
		\end{minipage}	
		\begin{minipage}[c][1\width]{0.48\textwidth}
			\hspace{5pt} 
			{\begin{overpic}[width=\textwidth]{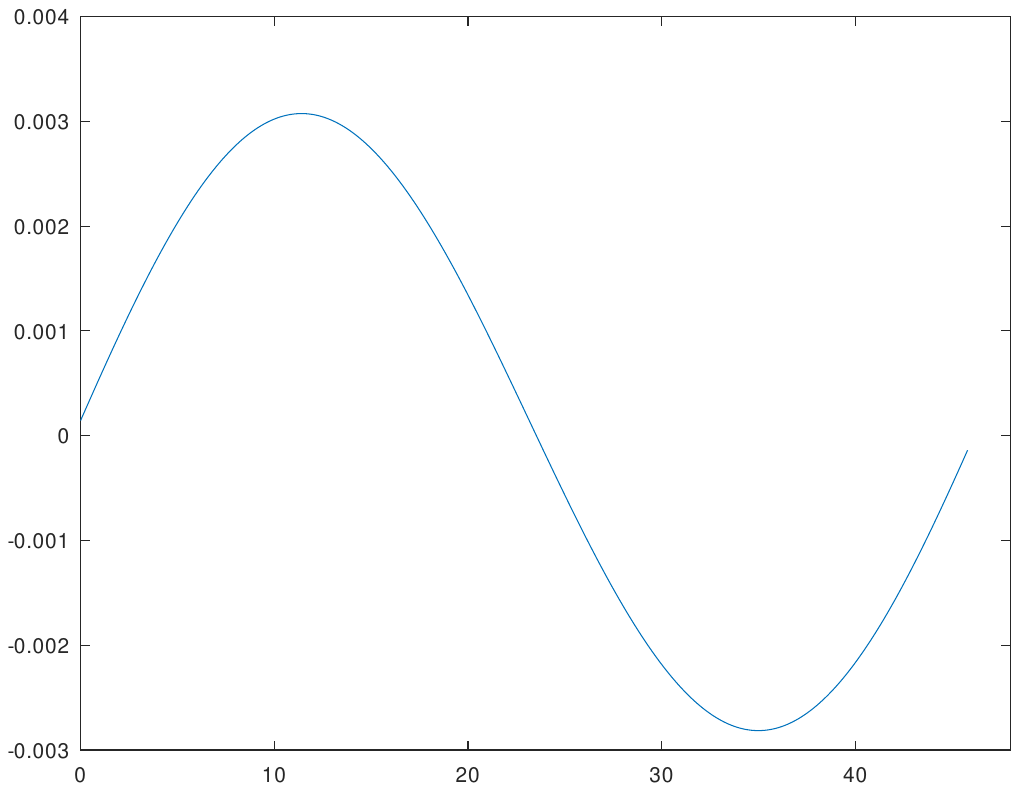}
					\put(0,78){\footnotesize (b)}
					\put(4,67.5){\footnotesize $\tilde{\beta}$}
					\put(63,24){\footnotesize $t$}
			\end{overpic}}
		\end{minipage}\\
		\vspace{-30pt}
  	\begin{minipage}[c][1\width]{0.48\textwidth}
			\hspace{5pt} 
			{\begin{overpic}[width=\textwidth]{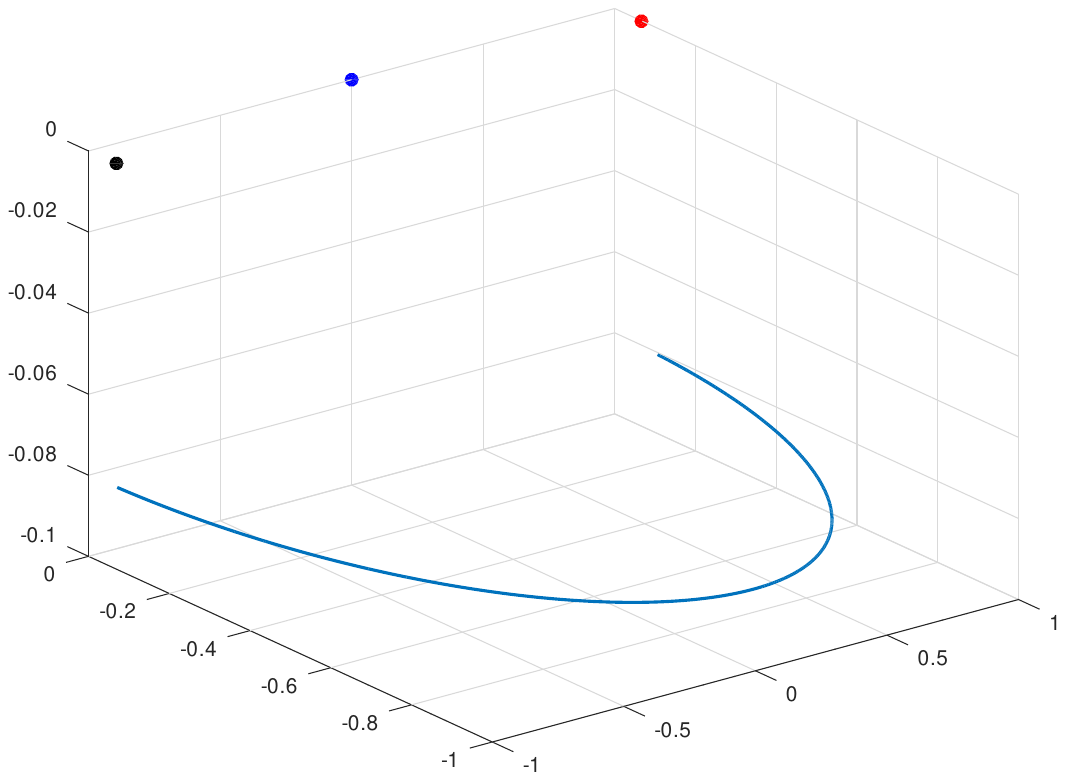}
					\put(0,78){\footnotesize (c)}
     	 		\put(13,67){\footnotesize $\ue{s}^{-}$}
         	        \put(27,71){\footnotesize $\ue{0}$}
                  	\put(44,75){\footnotesize $\ue{s}^{+}$}
         	     	\put(2,58){\footnotesize $u_3$}
     	            \put(17,28){\footnotesize $u_2$}
					\put(57,28){\footnotesize $u_1$}
			\end{overpic}}
		\end{minipage}	
    \begin{minipage}[c][1\width]{0.48\textwidth}
			\hspace{0pt}
			{\begin{overpic}[width=\textwidth]{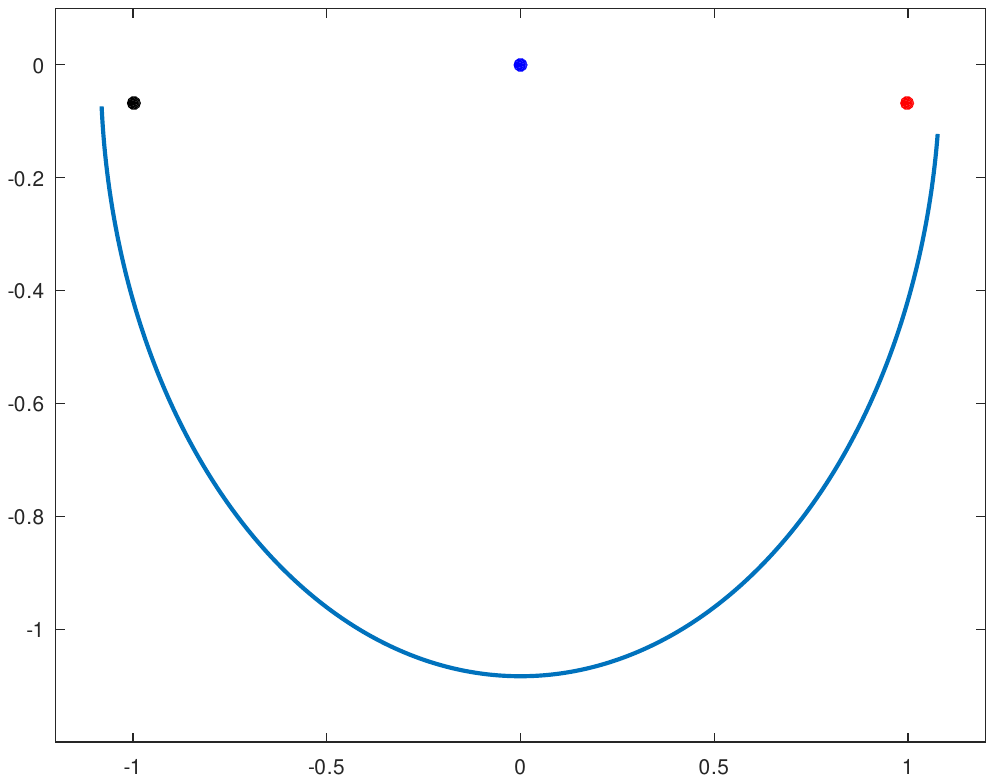}
					\put(0,78){\footnotesize (d)}
					\put(5,66){\footnotesize $w_1$}
					\put(57,24){\footnotesize $w_2$}
                        \put(15,69){\footnotesize $\ue{w}^-$}
                        \put(41,69){\footnotesize $\ue{0}$}
                        \put(62,69){\footnotesize $\ue{w}^+$}
 			\end{overpic}}
		\end{minipage}\\
	\end{center}
	\vspace{-30pt}
	\caption{A case study to show the quasi-latitudinal features of the constructed solutions during the ``transfer'' stage (the ``expulsuion'' stage is not reported for the sake of clarity). In panel (a) the behaviour of $u_3(t)$ for the values $\lambda_{(1)}=0.0055$, $\lambda_{(2)}=0.005$ and $\lambda_{(3)}=0.001$. The remaining panels show the case with $\lambda=\lambda_{(3)}$. Panel (b) shows the current realising the control, whilst panels (c) and (d) represent the trajectory of the system controlled with such a current in the space $\ue{u}$ and stereographic projection $\ue{w}$, respectively. In particular, $\ue{w}^{\pm}$ represent the stereographic projection of $\ue{s}^{\pm}$. The remaining parameters are reported in Appendix C.}
	\label{fig:transfer}
	\vspace{0pt}
\end{figure*}

Let us firstly introduce the following notation. Given any vector $\ue{x} \in \RR^3$ we shall denote with $\ue{x}_r:=(x_1,x_2)$. Vice-versa, we shall denote with $w_{1,2}$ the first two components of either $\ue{w}_r$ or $\ue{w}$. It is immediate to notice from (\ref{eq:remaindersone}) and (\ref{eq:remaindersonethree}) that  $\ue{\mathcal{F}}(\ue{v};\ue{v})\equiv (\ue{\mathcal{F}}_r(\ue{v}_r;\ue{v}_r),0)$ and similarly for $\ue{\ml{R}}$. Hence, it is meaningful to consider the following system 
\beq{eq:prswscaledthree}
\dot{\ue{w}}_r=\uue{L} \ue{w}_r+\lambda \ue{\ml{F}}_r(\ue{w}_r,\ue{w}_r) 
\eeq
where 
\beq{eq:sigmaomega}
\ue{L}:=\begin{pmatrix}
	0 & -\sigma \omega\\
	\omega / \sigma & 0\\
\end{pmatrix}
,\qquad \sigma:=\sqrt{D_{3,2}/D_{3,1}}\mx{,}
\eeq
and $\omega$ has been defined in (\ref{eq:omega}). Clearly, (\ref{eq:prswscaledthree}) is nothing but the first order truncation of (\ref{eq:prswscaledtwo}), written in the reduced set of variables $\ue{w}_r$. Moreover, it is immediate to check that 
\beq{eq:firstint}
\ml{G}_{tr}:=\sigma^{-1} w_1^2+\sigma w_2^2 \mx{,}
\eeq
is a prime integral for (\ref{eq:prswscaledthree}) if $\lambda=0$. \\
With the aim to construct CQL solutions, we ask ourselves if the system can be solved, for instance, by means of perturbative tools. The answer is affirmative, as stated in the following
\begin{prop}\label{prop:one}%
It is possible to construct a function $\ue{w}_r^{[\leq 1]}(t)$ satisfying (\ref{eq:prswscaledthree}) up to $O(\lambda)$. 
\end{prop}
\begin{rem}
It is important to avoid any ambiguity about the meaning of the previous statement. Solving (\ref{eq:prswscaledthree}) up to $O(\lambda)$ does not mean solving it exactly. In fact one has 
$\dot{\ue{w}}_r^{[\leq 1]}-\uue{L} \ue{w}_r^{[\leq 1]}-\lambda \ue{\ml{F}}_r(\ue{w}_r^{[\leq 1]},\ue{w}_r^{[\leq 1]})=O(\lambda^2)$. This is typical of perturbative arguments, which are known to generate a remainder as a consequence of the expansions involved.
\end{rem}
\proof
Let us firstly cast the linear part of (\ref{eq:prswscaledthree}) into a diagonal form. It is immediate to check that the required transformation is given by
\beq{eq:coordinatestransf}
\ue{w}_r=\ue{C} \ue{x},\qquad
\ue{C}:=\begin{pmatrix}
	\sigma & \sigma\\
	-i & i\\
\end{pmatrix} \mx{.}
\eeq
In fact, $\ue{C}^{-1} \ue{L} \ue{C}=\diag (i \omega, -i \omega)=:\ue{\Lambda}$. We can now cast system  (\ref{eq:prswscaledthree}) in the new set of variables
\beq{eq:systemnewset}
\dot{\ue{x}}=\ue{\Lambda} \ue{x}+ \lambda \ue{C}^{-1}  \ue{\mathcal{F}}_r (\ue{C} \ue{x}) \mx{,} 
\eeq
into a Hamiltonian form via a phase space extension. More precisely, by denoting with $y_{j}$ the momenta canonically conjugated to $x_j$, one has that (\ref{eq:systemnewset}) is given by (part of) the canonical equations of 
\[
H(\ue{y},\ue{x})=H_0(\ue{y},\ue{x})+\lambda H_1(\ue{y},\ue{x}) \mx{,}
\]
where
\beq{eq:hone}
\begin{aligned}
H_0 & := i \omega (x_1 y_1 - x_2 y_2) \\
H_1 & :=K \rho y_1  \left[\frac{\thd}{2} \left(\sigma+\frac{1}{\sigma} \right) x_1^2 + \frac{\thd}{2} \left(\frac{1}{\sigma}-\sigma \right) x_2^2 
-  \frac{\thd}{\sigma} x_1 x_2 - i \tilde{D}_{2,1} \sigma x_1 x_2^2 + i \tilde{D}_{2,1} \sigma x_1^3 
\right]\\
& + K \rho y_2  \left[\frac{\thd}{2} \left(\sigma+\frac{1}{\sigma} \right) x_2^2 + \frac{\thd}{2} \left(\frac{1}{\sigma}-\sigma \right) x_1^2 
-  \frac{\thd}{\sigma} x_1 x_2 + i \tilde{D}_{2,1} \sigma x_1^2 x_2 - i \tilde{D}_{2,1} \sigma x_2^3 
\right] 
\end{aligned} \mx{.}
\eeq 
As it is common in perturbation theory, we ask weather it is possible to find a canonical, $\lambda-$close to the identity transformation of variables $(\ue{x},\ue{y}) = \mathcal{T} (\ue{X},\ue{Y})$ apt to ``remove'' the contribution of $H_1$ i.e. such that 
\beq{eq:normalformcast}
H \circ \mathcal{T}=H_0(\ue{Y},\ue{X})+O(\lambda^2) \mx{.}
\eeq
For this purpose, by invoking the well known Gr\"{o}bner Exchange Theorem, see e.g. \cite{gior03}, such a transformation will be determined by requiring that
\[
\exp( \mathcal{L}_{\lambda \chi}) H = H_0+O(\lambda^2) \mx{,}
\]
i.e. the well known first order homological equation
\beq{eq:homological}
\{ \chi,H_0\}=H_1 \mx{,}
\eeq
where the generating function $\chi$ will be sought of the form  
\beq{eq:genfun}
\chi(\ue{y},\ue{x})=\ue{y} \cdot \ue{\mathcal{C}} (\ue{x}) \mx{.}
\eeq
This choice is suggested by a general property when dealing with Hamiltonians obtained via the above described phase space extension, see \cite{for162} for a proof. In this particular case, we will look for $\mathcal{C}_j$ as non-homogeneous polynomials of the form
\beq{eq:formc}
\mathcal{C}_j(\ue{x}):=\sum_{|\ue{\nu}|=2,3}c_{\ue{\nu}}^{(j)}\ue{x}^{\ue{\nu}} \mx{,}
\eeq
where $\ue{\nu} \in \NN^2$, $|\nu|:=\nu_1+\nu_2$, $\ue{x}^{\ue{\nu}}:=x_1^{\nu_1} x_2^{\nu_2}$ and $c_{\ue{\nu}}^{(j)}$ are complex-valued unknown coefficients to be determined. By using (\ref{eq:formc}) in (\ref{eq:genfun}), then substituting in (\ref{eq:homological}) one finds
\beq{eq:coeffdegreetwoa}
\begin{aligned}
-c_{0,2}^{(1)}=c_{2,0}^{(2)}&=i \thd K \rho (\sigma^2-1)/(6 \sigma \omega ) \\ 
-c_{2,0}^{(1)}=c_{0,2}^{(2)}&=i \thd K \rho (\sigma^2+1)/(2 \sigma \omega ) \\
-c_{1,1}^{(1)}=c_{1,1}^{(2)}&=i \thd K \rho \sigma / (\sigma \omega)
\end{aligned}
\eeq
and 
\beq{eq:coeffdegreetwob}
\begin{aligned}
c_{3,0}^{(1)}=c_{0,3}^{(2)}=c_{1,2}^{(1)}=c_{2,1}^{(1)}&=\tilde{D}_{2,1}K \rho \sigma /(2 \omega )\\
c_{3,0}^{(1)}=c_{0,3}^{(2)}&=0
\end{aligned}
\eeq
The clear symmetry relations amongst these coefficients is related to the well known property $x_1=\bar{x}_2$ (here $\bar{z}$ denotes the complex-conjugated of $z \in \CC$), which is typical of the coordinates maps (\ref{eq:coordinatestransf}). In conclusion, if $O(\lambda^2)$ terms are disregarded, the normalising transformation reads as
\beq{eq:transft}
\ue{x}=\ue{\mathcal{N}}_{[\leq 1]}(\ue{X}):=\ue{X}+\lambda \ue{\mathcal{C}}(\ue{X}) \mx{.}
\eeq
If the above mentioned Gr\"{o}bner Theorem is not used, checking that (\ref{eq:transft}) satisfies (\ref{eq:normalformcast}) is just a matter of patience. 
\begin{rem}
Interestingly enough, the variables $\ue{y}$ and $\ue{Y}$ do not appear in (\ref{eq:transft}), which consistently becomes a normalising transformation in the pair $(\ue{x},\ue{X})$ alone for the original system (\ref{eq:systemnewset}).
\end{rem}
It is a simple consequence of the construction above that if $\ue{X}(t)$ is a solution to the normalised system (up to $O(\lambda)$), then $\ue{x}(t):=\ue{\mathcal{N}}_{[\leq 1]}(\ue{X}(t))$ is a solution to the original system (up to $O(\lambda)$).\\
On the other hand, the integration of the normalised system is immediate
\[
X_{1,2}(t)=(X_a \pm i X_b)e^{\pm i t \omega} \mx{,}
\] 
where 
\beq{eq:xab}
X_a:=\Re X_1(0), \qquad X_b:=\Im X_1(0) \mx{,} 
\eeq
(recall that $X_1=\bar{X}_2$). Hence, by using (\ref{eq:coordinatestransf}), one gets
\beq{eq:integrationv}
\begin{aligned}
w_{1}(t) &= 2 \sigma (X_a \cos(t \omega)-X_b \sin(t \omega)) + \lambda \omega^{-1} K \rho\left\{ \tilde{D}_{2,1}  \sigma^2 (X_a^2+X_b^2) [X_a \cos(t \omega) - X_b \sin(t \omega)] \right.  \\
&+(2/3) \thd (\sigma^2+2)[2 X_a X_b \cos(2 t \omega)-(X_b^2-X_a^2) \sin(2 t \omega) ] \\
& \left.- \tilde{D}_{2,1} \sigma^2 [X_a (3 X_b^2-X_a^2)\cos(3 t \omega) - X_b (X_b^2- 3 X_a^2)\sin(3 t \omega) ] \right\}    \\
w_{2}(t) &= 2 (X_b \cos(t \omega)+X_a \sin(t \omega)) + \lambda \omega^{-1} K \rho\left\{ - \tilde{D}_{2,1}  \sigma (X_a^2+X_b^2) [X_b \cos(t \omega) + X_a \sin(t \omega)] \right.  \\
&+(2/3) \thd (2 \sigma^2+1)[(X_b^2-X_a^2) \cos(2 t \omega)+ 2 X_a X_b \sin(2 t \omega) ] \\
& \left.- \tilde{D}_{2,1} \sigma [X_b (X_b^2- 3 X_a^2) \cos(3 t \omega) + X_a  (3 X_b^2-X_a^2)\sin(3 t \omega) ] - 2 \thd (X_a^2+X_b^2) / \sigma \right\}    \\
\end{aligned} \mx{.}
\eeq
Clearly, the functions above provide the required approximated solution to a given Cauchy problem for (\ref{eq:prswscaledthree}), once the quantities $X_{a,b}$ have been determined from the initial condition $\ue{w}_r(0)$ (which will be specified later on). More precisely, by (\ref{eq:coordinatestransf}) and (\ref{eq:transft}) one has
\beq{eq:xzero}
\ue{X}(0)=\ue{\mathcal{N}}_{[\leq 1]}^{-1}(\ue{C}^{-1} \ue{w}_r(0)) \mx{,}
\eeq
where we recall that, if $O(\lambda^2)$ terms are disregarded, one simply has $\ue{\mathcal{N}}_{[\leq 1]}^{-1}(\ue{x})=\ue{x}-\lambda \ue{\mathcal{C}}(\ue{x})$. The values of $X_{a,b}$ are finally determined by recalling (\ref{eq:xab}). These quantities determine the required $\ue{w}_r^{[\leq 1]}(t)$ via (\ref{eq:integrationv}).   
\endproof
\begin{rem}
Clearly, the normalisation order could be increased by generalising the above described procedure, for instance by using the Lie Series or Lie Transform methods, see e.g. \cite{gior03}. The Hamiltonian formulation plays a key role from this point of view.\\
Knowingly, the computations get dramatically more and more involved as the normalisation order is increased. 
\end{rem}
Let us now notice that at the zero-th order in $\lambda$, (\ref{eq:xzero}) yields 
\[
\ue{X}^{[0]}(0)=2^{-1}(\sigma^{-1} w_{1}(0),w_{2}(0)) \mx{,}
\]
hence, from (\ref{eq:integrationv}), one gets
\[
\ue{w}_r^{[0]}(t)=(w_{1}(0)\cos (\omega t)- w_{2}(0) \sigma \sin (\omega t),w_{1}(0) \sigma^{-1} \sin (\omega t)+w_{2}(0)\cos (\omega t)) \mx{.}
\]
As $w_{1}(0)$ will be assumed as negative and bounded away from zero, the first component can be written as $\ml{A}_{m}\cos(\omega t+\phi)$, where 
$\phi:=\arctan (\sigma w_{2}(0)/w_1(0))$ and 
\beq{eq:am}
\ml{A}_m:=w_{1}(0)\sqrt{1+(\sigma w_{2}(0)/w_1(0))^2}<0\mx{.} 
\eeq
Hence, we shall define $T_{tr}$ as the time necessary for $w_{r,1}^{[0]}(t)$ to reach $-\ml{A}_m-K^2$. This yields 
\beq{eq:ttr}
T_{tr}:=\omega^{-1} \left[\arccos(-1-K^2/\ml{A}_m)-\arctan(\sigma w_{2}(0)/w_1(0)) \right] \mx{.}
\eeq
It is immediate to realise that this quantity is correctly defined as we we shall only consider cases in which $K^2 \ll |w_1(0)|\leq |\ml{A}_m|$. Moreover, recalling (\ref{eq:ttrplus}), it is easy to check that  
\beq{eq:upperttr}
T_{tr} \leq T_{tr}^+ \mx{.}
\eeq
In fact, $T_{tr}^+$ is attained in the (unrealistic as $D_{2,1}>0$) case $\ue{w}_r(0)=(-1,0)$ i.e $\gamma=1$. However, as it is evident from (\ref{eq:translation}) and (\ref{eq:solutiontoexp}), the expulsion stage acts by increasing $w_{2}(t)$ (and hence decreasing $w_1(t)$, and this clearly reduces $T_{tr}$.\\
Once the control has been constructed, it is important to obtain an estimate of the measure of the initial conditions around $\ue{w}(0) \equiv (\ue{w}_r(0),-K)$ which are ``safely transported'' in a neighbourhood of the target point $\ue{w}(T_{tr})$. In order to achieve this, we need to bound the difference between the solutions of two different systems. The first one is given by (\ref{eq:prswscaledthree}) (i.e. the $O(\lambda)$ truncation of (\ref{eq:prswscaledtwo})), whose solution is known by Prop. \ref{prop:one}. By setting $\ue{w}:=(\ue{w}_r^{[\leq 1]},-K)$ and disregarding $O(\lambda^2)$, this reads as 
\beq{eq:systemw}
\dot{\ue{w}}=:\uue{A}(-K)\ue{w}+\lambda \ue{\mathcal{F}}(\ue{w};\ue{w}) \mx{,}
\eeq
being subject to the prefixed initial condition $\ue{w}(0)$. The second system is the full, controlled one, i.e. obtained from (\ref{eq:prswscaled}) by setting $\tilde{\beta}(t):=\hat{\beta}_{lat}(\ue{w}^{[\leq 1]}(t))$. Its $O(\lambda)$ truncation can be written as
\beq{eq:newsys}
\dot{\ue{u}}=\uue{A}(u_3)\ue{u}+\lambda \ue{\mathcal{F}}(\ue{u};\ue{w})  \mx{,}
\eeq 
with an ``uncertain'' initial condition around $\ue{w}(0)$. Let us denote by $\Phi_{tr}^{t}$ the flow of the controlled system (\ref{eq:newsys}). An example of the constructed CQL solutions and the corresponding control is shown in fig. \ref{fig:transfer}.   
\begin{lem}[$O(\lambda)$ cut-off]\label{lem:transfer}
Let us define the following quantities
\[
\ml{M}_{tr} :=1+\max_{|\ue{u}|\leq 2 } ||\uue{D_{\ue{u}} \ml{F}}(\ue{u};\ue{w})||_{\infty}, \qquad    
\ml{K}_{\ue{w}}:=1+\max_{t \in [0,T_{tr}]} |\ue{w}_r^{[\leq 1]}| \mx{.}
\]
Then, for all $\lambda \in (0,\lambda_{tr}]$, with 
\beq{eq:hypk}
\lambda_{tr} := (4 \pi^2 \ml{K}_{\ue{w}} \ml{M}_{tr} )^{-1}\omega^2 \mx{,}
\eeq
and all $\rho_{tr}^+>0$, there exists $\rho_{tr}^-=\rho_{tr}^-(\rho_{tr}^+)$, such that following property holds
\beq{eq:inclusion}
\Phi_{tr}^{T_{tr}}(\mathfrak{B}_{\lambda \rho_{tr}^-}(\ue{w}(0)) \subset
\mathfrak{B}_{\lambda \rho_{tr}^+}(\ue{w}(T_{tr})) \mx{.}
\eeq
\end{lem} 
\begin{rem}\label{rem:threshold}
The $O(\lambda^2)$ terms in (\ref{eq:newsys}) can be bounded by using the Lie series theory, and this requires some extra work. However, it is possible to check that the information carried by them can be thought as ``negligible''. More precisely, by facing slightly more cumbersome estimates, one could repeat the proof below by including these $O(\lambda^2)$ contributions, and verify that those terms would imply an additional $O(\lambda)$ term in (\ref{eq:hypk}). Hence, the above statement still holds for the ``non cut-off'' problem at a price of a possible (but ``small'') further restriction of the threshold $\lambda_{tr}$.  
\end{rem}
\proof This proof has some similarities with the one of Lem. \ref{lem:expulsion}. Let us now define $\ue{\eta}:=\ue{u}-\ue{w}$, hence, by (\ref{eq:systemw}) and (\ref{eq:newsys}) one gets
\beq{eq:etaode}
\dot{\ue{\eta}}=\uue{A}(-K)\ue{\eta}+\uue{A}(\eta_3)(\ue{w}+\ue{\eta})+
\lambda[\ue{\ml{F}}(\ue{w}+\ue{\eta};\ue{w})-\ue{\ml{F}}(\ue{w};\ue{w})] \mx{.}
\eeq
It is evident from the latter and the structure of $\ue{A}$ that $\dot{\eta}_3$ has an $O(\lambda)$ magnitude, hence it is appropriate to consider the new variable $\eta_3=:\lambda \zeta_3$. By substituting it in (\ref{eq:etaode}), one gets that the $O(1)$ evolution of $\eta_{1,2}$ is given by the linear system $\dot{\ue{\eta}}_{r}=\uue{L} \ue{\eta}_{r}$, where $\ue{\eta}_{r}:=(\eta_1,\eta_2)$. \\
Let us now introduce $d_{r}^{\pm},d_3^{\pm} \in (0,1]$, which will play the role of real parameters to be determined, then suppose for a moment that the following bound holds
\beq{eq:ansatzzeta}
|\zeta_3| \leq d_3^+ \mx{.}
\eeq
By multiplying both sides of the first two components of equation (\ref{eq:etaode}) by the ``integrating factor'' $\exp(-\uue{L}t)$, integrating and finally considering the absolute values, one obtains
\[
\begin{aligned}
|\ue{\eta}_r(t)| & \leq  |\ue{\eta}_r(0)| + \lambda \int_0^t 
\left\{|\ue{A}(\zeta_3)(\ue{w}+\ue{\eta})|+ \ml{M}_{tr} |\ue{\eta}| \right\} ds \\
& \leq |\ue{\eta}_r(0)| + \lambda d_3^+ \ml{K}_{\ue{w}} t +\lambda  (d_3 + \ml{M}_{tr})  \int_0^t |\ue{\eta}_r(s)|ds
\end{aligned}
\]
where we have used that $|\exp(\uue{L}\cdot) \ue{z}|=|\ue{z}|$, for all $\ue{z} \in \RR^3 $ as $\uue{L}$ has pure complex eigenvalues only, and the bound
\beq{eq:boundf}
|\ml{F}_i(\ue{w}+\ue{\eta};\ue{w})-\ml{F}_i(\ue{w};\ue{w})|\leq|\ue{\ml{F}}(\ue{w}+\ue{\eta};\ue{w})-\ue{\ml{F}}(\ue{w};\ue{w})| \leq \ml{M}_{tr} |\ue{\eta}| \mx{.}
\eeq 
Hence, by using the Gr\"{o}nwall lemma,
\beq{eq:etatwodim}
|\ue{\eta}_r(t)| \leq \left[  |\ue{\eta}_r(0)| + \lambda d_3^+ \ml{K}_{\ue{w}} t \right] 
e^{\lambda  (d_3^+ + \ml{M}_{tr}) t} \mx{.}
\eeq
Recalling that $T_{tr} \leq \pi/\omega$ by (\ref{eq:upperttr}), by using assumption (\ref{eq:hypk}) one gets $\lambda  (1 + \ml{M}_{tr}) t \leq  2 \lambda  \ml{M}_{tr} t  \leq \log(2)$ for all $t \leq T_{tr}$ as $\ml{M}_{tr} \geq 1$. Hence, in order to obtain 
\beq{eq:drp}
|\ue{\eta}_r(t)| \leq \lambda d_r^+ \mx{,}
\eeq
it is sufficient to require $|\ue{\eta}_r(0)| \leq \lambda d_r^-$, and then
\beq{eq:conddtr}
d_r^- \leq d_r^+/4,\qquad d_3^+ \leq \omega d_r^+/(\pi \ml{K}_{\ue{w}}) \mx{.}
\eeq
Let us now consider the third equation of (\ref{eq:etaode}), i.e. $\dot{\zeta}_3=\ml{F}_3(\ue{w}+\ue{\eta};\ue{w})-\ml{F}_3(\ue{w};\ue{w})$. The latter and bound (\ref{eq:boundf}), yield $|\zeta_3(t)|\leq |\zeta_3(0)|+\int_0^t \ml{M}_{tr} \left(|\ue{\eta}_r(s)|+\lambda |\zeta_3(s)| \right)ds $. Hence, (\ref{eq:drp}) yields 
\[
|\zeta_3(t)| \leq |\zeta_3(0)| + \lambda 
\ml{M}_{tr} d_r^+ t + \lambda \ml{M}_{tr} \int_0^t |\zeta_3(s)|ds \mx{.}
\]
By using the Gr\"{o}nwall lemma once again and assumption (\ref{eq:hypk}), which implies, \emph{a fortiori}, $\exp(\lambda \ml{M}_{tr} t)<2$ for all $t \in [0, T_{tr}]$, we obtain 
\[
|\zeta_3(t)| \leq 2 \left[ |\zeta_3(0)| + \lambda \ml{M}_{tr} d_r^+ \pi/\omega \right] 
\mx{.}
\]
Hence, by setting $|\zeta_3(0)| \leq d_3^-$, a sufficient condition for the property (\ref{eq:ansatzzeta}) to hold true, is  
\beq{eq:conddtrtwo}
d_3^-\leq d_3^+/4, \qquad \lambda d_r^+ \leq \omega d_3^+/(4 \pi \ml{M}_{tr} ) \mx{.}
\eeq
It is easy to realise that the second conditions appearing in (\ref{eq:conddtr}) and in (\ref{eq:conddtrtwo}) hold in the region of the plane $(d_r^+,d_3^+)$ given by $d_3^+ \in [4 \omega^{-1}\pi \ml{M}_{tr} \lambda d_r^+, \pi^{-1} K_{\ue{w}}^{-1} \omega d_r^+]$, for all $d_r^+\in[0,1]$. Such a region is correctly defined by virtue of (\ref{eq:hypk}). Hence, for any $d_r^+$, we choose the medium point of the above mentioned interval, that is, 
\beq{eq:drtp}
d_3^+:=\Theta d_r^+ \mx{,}
\eeq
where $\Theta:=2^{-1} (\omega^{-1} \lambda 4 \pi \ml{M}_{tr}+ \pi^{-1} K_{\ue{w}}^{-1} \omega )$. Note that $\Theta<1$ (as, in particular, $(2\pi)^{-1} K_{\ue{w}}^{-1} \omega<1$ by construction). This implies $d_{3}^+<d_r^+$. Hence, as $|\ue{\eta}(0)|\geq \max\{|\ue{\eta}_r(0)|,|\eta_3(0)|\}$, the first conditions of (\ref{eq:conddtr}) and (\ref{eq:conddtrtwo}) hold if one requires
\beq{eq:lhstr}
\lambda(d_3^+/4) \geq |\ue{\eta}(0)| \mx{.}  
\eeq
On the other hand, bounds (\ref{eq:ansatzzeta}) and (\ref{eq:drp}) imply, for all $t \leq T_{tr}$, 
\beq{eq:lasttr}
|\ue{\eta}(t)| \leq |\ue{\eta}_r(t)|+\lambda |\zeta_3(t)| \leq \lambda \left[1+\Theta \right]  d_r^+ \mx{.}
\eeq
In conclusion, given $\rho_{tr}^+>0$, we have that (\ref{eq:lasttr}) and (\ref{eq:lhstr}), imply that by defining   
\beq{eq:choicerhotr}
\rho_{tr}^-:= [4(1+\Theta)]^{-1} \Theta \rho_{tr}^+ \mx{,}
\eeq
and choosing $|\ue{\eta}(0)| \leq 
\lambda \rho_{tr}^-$ as in (\ref{eq:inclusion}), the proof is complete. 
\endproof

\section{Attraction}\label{sec:five}
\begin{figure*}[t!]\begin{center}
  		\begin{minipage}[c][1\width]{0.48\textwidth}
			\hspace{0pt}
			{\begin{overpic}[width=1\textwidth]{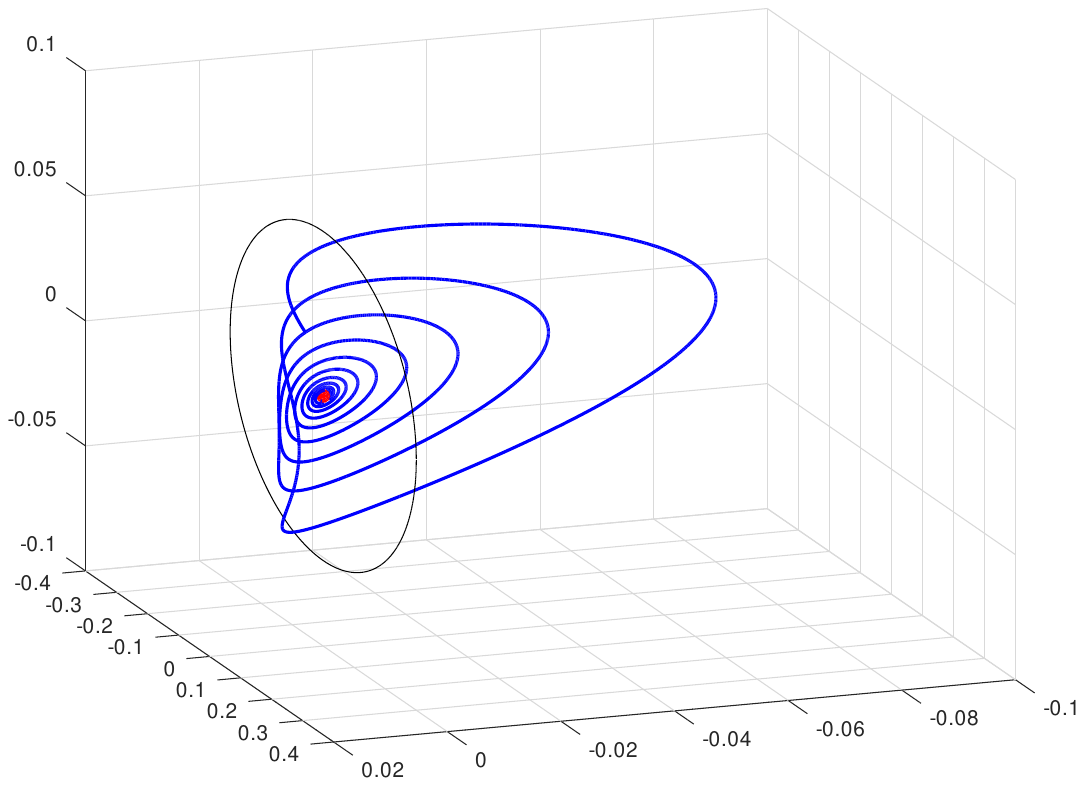}
					\put(0,78){\footnotesize (a)}
        	     	\put(3,59){\footnotesize $U_3$}
     	            \put(10,30){\footnotesize $U_2$}
					\put(60,26){\footnotesize $U_1$}
			\end{overpic}}
		\end{minipage}	
		\begin{minipage}[c][1\width]{0.48\textwidth}
			\hspace{5pt} 
			{\begin{overpic}[width=\textwidth]{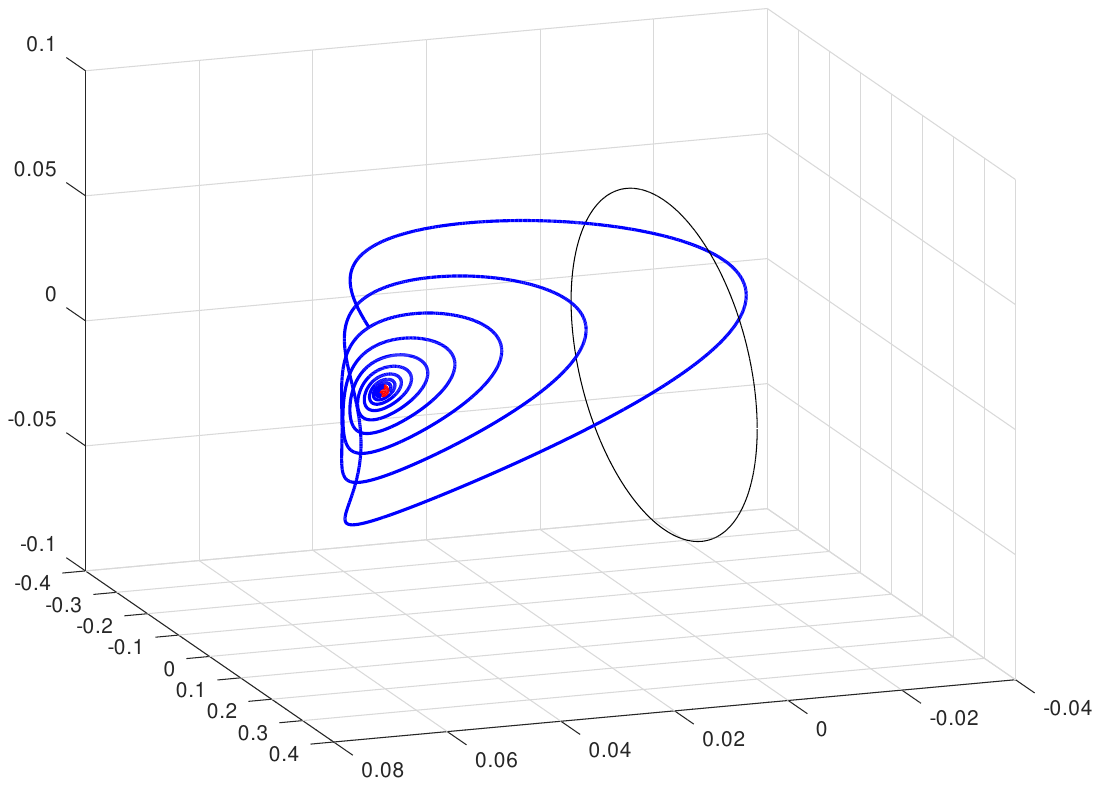}
					\put(0,78){\footnotesize (b)}
        	     	\put(3,59){\footnotesize $U_3$}
     	            \put(10,30){\footnotesize $U_2$}
					\put(60,26){\footnotesize $U_1$}
			\end{overpic}}
		\end{minipage}\\
		\vspace{-30pt}
  	\begin{minipage}[c][1\width]{0.48\textwidth}
			\hspace{5pt} 
			{\begin{overpic}[width=\textwidth]{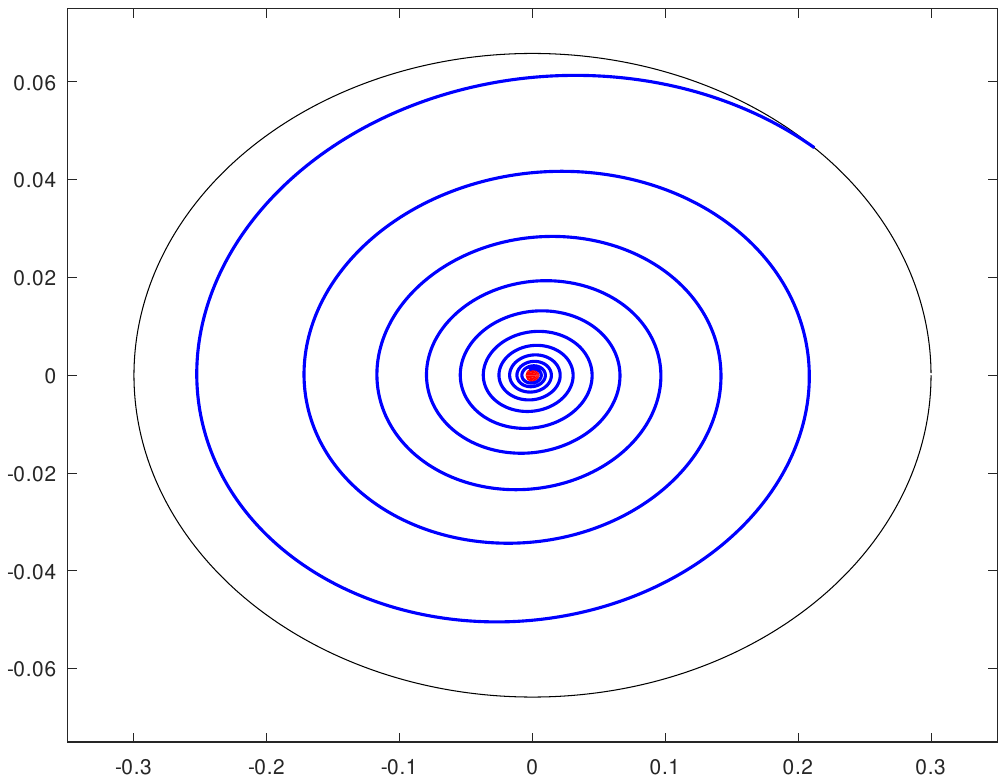}
					\put(0,78){\footnotesize (c)}
         	     	\put(3,60){\footnotesize $U_3$}
     	            \put(51,24){\footnotesize $U_2$}
			\end{overpic}}
		\end{minipage}	
    \begin{minipage}[c][1\width]{0.48\textwidth}
			\hspace{0pt}
			{\begin{overpic}[width=\textwidth]{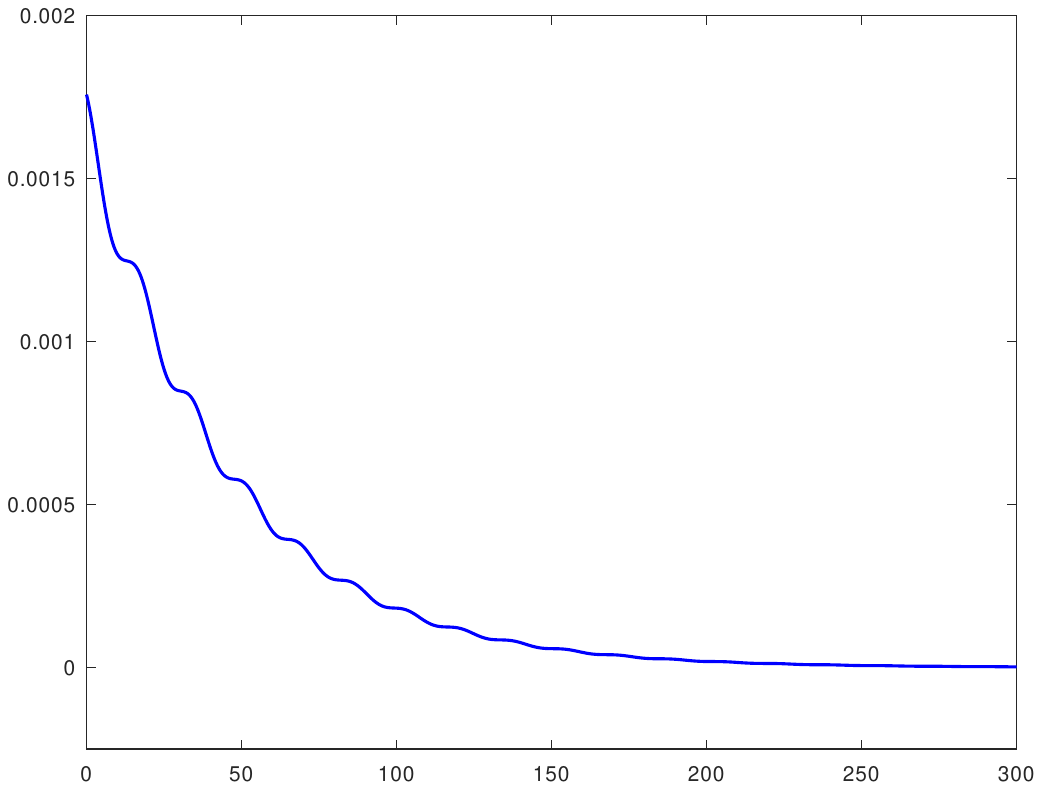}
					\put(0,78){\footnotesize (d)}
					\put(3,67){\footnotesize $\mathfrak{W}$}
					\put(53,24){\footnotesize $t$}
			\end{overpic}}
		\end{minipage}\\
	\end{center}
	\vspace{-30pt}
	\caption{The field attractivity feature within the basin $\ml{B}_{\ml{E}}$ when the injected current $\beta(t)$ is switched off. The (black) ellipse in the first three panels represents the set $\mathfrak{W}(U_2,U_3) = (D_{2,1} \gamma)^2/(32 D_{3,1} \Omega^2)$ (recall that $\ue{U}$, defined via (\ref{eq:translsplus}), represents a coordinate system in which the equilibrium $\ue{s}^+$ corresponds to the origin). Panels (a) and (b) show a trajectory starting on the boundary of $\ml{B}_{\ml{E}}$ for $\delta_a=0$ and $\delta_a=0.1$, respectively. In panel (c) the projection on the $(U_2,U_3)$ plane of the latter trajectory is depicted. The last panel shows, as expected, the decreasing behaviour of the Lyapunov function $\mathfrak{W}$ along the same trajectory. The full set of parameters is reported in Appendix C.}
	\label{fig:attraction}
	\vspace{0pt}
\end{figure*}

As anticipated, the aim of this section is to study the dynamics in a neighbourhood of the equilibrium $\ue{s}^+$ when the current is switched off, i.e. $\beta=0$, and show that $\ue{s}^+$ is an attractive point for the system, providing in addition an estimate for the basin of attraction. For this purpose, we shall proceed in a non-perturbative fashion and consider system (\ref{eq:prswscaledalternative}). This is motivated by the fact that the attractive behaviour is characterised by $\alpha$, whilst $\lambda$ does not play a particularly relevant role in this case. As usual, we shall start by considering the standard translation  
\beq{eq:translsplus}
\ue{u}=\ue{U}+\ue{s}^+ \mx{,}
\eeq  
and setting $\beta=0$ then $h_2=:-D_{2,1}\Omega$ with $\Omega>0$, and $\gamma:=\sqrt{1-\Omega^2}$, in such a way (\ref{eq:prsw}) reads as
\beq{eq:systemu}
\dot{\ue{U}}=\ue{\ml{G}}_1(\ue{U})+\alpha \ue{\ml{G}}_2 (\ue{U}) \mx{,}
\eeq
where $\ue{\ml{G}}_1(\ue{U}):=(D_{3,1} U_3 (U_2-\Omega) - D_{3,2} U_2 U_3 ,-D_{3,1} U_3 (\gamma + U_1), D_{2,1} U_2 (\gamma+ U_1)   )^{\top}$ and $\ue{\ml{G}}_2$ is defined as a consequence. We shall denote with $\Phi_a^t$ the corresponding phase flow. Moreover, let us define $\tilde{\Psi}(\ue{U}):=\Psi(\ue{u})|_{\ue{u}=\ue{U}+\ue{s}^+}$. \\ 
In this setting, we can state the following
\begin{lem}\label{lem:attraction} 
Define the following quadratic function 
\beq{eq:liapf}
\mathfrak{W}(U_2,U_3):=2^{-1}(D_{2,1}U_2^2+D_{3,1}U_3^2) 
\eeq
and suppose 
\beq{eq:choiceomegak}
3 \Omega^2 \geq  2 D_{2,1}/D_{3,1}, \qquad 16 \sqrt{D_{2,1}/D_{3,1}} \leq \gamma\mx{.}
\eeq
Then, for all $\alpha>0$, and all $\ml{E}:=1+\delta_a$, with 
\beq{eq:conddeltaa}
|\delta_a|\leq (\gamma/4)^2\mx{,}
\eeq
all the solutions of (\ref{eq:systemu}) starting in the set  
\[
\ml{B}_{\mathcal{E}}:=\{\ue{U}\in \RR^3:|U_1|\leq \gamma/4, \quad \mathfrak{W}(U_2,U_3) \leq (D_{2,1} \gamma)^2/(32 D_{3,1} \Omega^2)\} \cap \{  \tilde{\Psi}(\ue{U})=\mathcal{E} \} \mx{,}
\]
satisfy the following property
\beq{eq:asymptotic}
\ue{U}_{\infty}(\delta_a):=\lim_{t \rw + \infty} \Phi_a^t \left( \ml{B}_{\mathcal{E}} \right) = (-\gamma +\sqrt{\gamma^2+\delta_a},0,0) \mx{,}
\eeq
and 
\beq{eq:absu}
|\ue{U}_{\infty}(\delta_a)| \leq \delta_a/(2 \gamma) \mx{.}
\eeq
In particular, those trajectories with  $\ml{E}=1$ are asymptotic to $\ue{U}=\ue{0}$ (and hence $\ue{u}(t) \rw \ue{s}^+$).\\
Furthermore, the projection of $\ml{B}_{\ml{E}}$ on the plane $(U_2,U_3)$ contains the disk centred at the origin of radius
\beq{eq:semiminor}
r_{sm}:=(4|\Omega|D_{3,1})^{-1}\gamma D_{2,1} \mx{.} 
\eeq
\end{lem} 
\begin{rem}
The importance of these last two lines of statement lies in the fact that, by construction, the CQL solution which realises the transfer, ends up at a point located in the vicinity of $\ue{U}=(0,0,-K)$. This will be made precise later (see, for instance, (\ref{eq:wte})). However, it is important to stress for the moment that this is necessary and we shall choose $2K \leq r_{sm}$ later on, in order to ensure that such a point is suitably attracted once the current is switched off.    
\end{rem}
\noindent A numerical validation of Lem. \ref{lem:attraction} is reported in fig. \ref{fig:attraction}.
\proof The proof relies on the Theory of Lyapunov functions and related tools for the basin of attraction estimation, see e.g. \cite{khal} for a comprehensive description. The proof will be carried out in the variables $U_2,U_3$, being the behaviour of the first variable a consequence, by the conservation law $\Psi(\ue{u})=\const$, see (\ref{eq:energy}).\\ 
Let us start noticing that (\ref{eq:liapf}) satisfies $\mathfrak{W}(0,0)=0$ and $\mathfrak{W}(U_2,U_3)>0$ for all $(U_2,U_3) \neq \ue{0}$. Furthermore, it is immediate to check that its derivative along the solutions of the ``undamped'' (\ref{eq:systemu}), i.e. with $\alpha=0$, satisfies
\[
\dot{\mathfrak{W}}|_{\alpha=0} := \nabla_{(U_2,U_3)}\mathfrak{W} \cdot \ue{\ml{G}}_1(\ue{U})=0 \mx{,}
\] 
hence $\mathfrak{W}$ is a (non-strict) Lyapunov function for the undamped (partial) system, showing that it is stable. However, our aim is to show that the damped system satisfies $\dot{\mathfrak{W}}  <0$ on a suitable set, proving in this way the (stronger) asymptotic stability.\\
As $D_{2,1} \ll D_{3,1}$ by assumption, the region we are considering to prove the negative definiteness of $\dot{\mathfrak{W}}$ will be conveniently chosen as non-isotropic. More precisely, we shall proceed by introducing in $\dot{\mathfrak{W}}$ the following variables transformation 
\[
\ue{U}=(\mu  x, \ep  (\gamma D_{2,1})^{-1} \cos \theta, \ep  D_{3,1}^{-1} \sin \theta) \mx{,}
\]
with $(x,\theta) \in [0,1] \times [0,2 \pi]$ and $\ep,\mu>0$ to be determined. The latter gives
\beq{eq:derivativew}
\begin{aligned}
\dot{\tilde{\mathfrak{W}}} &=-\alpha \ep^2\left[1+ 2 \mu x(\sin^2 \theta + \gamma^{-1} \cos^2 \theta)+ (\mu x)^2(\sin^2 \theta + \gamma^{-2} \cos^2 \theta) + \right. \\ 
& \left. - 2 \Omega \ep \gamma^{-1} (D_{2,1}^{-1}-D_{3,1}^{-1})\cos \theta \sin^2 \theta 
+ \ep^2 \gamma^{-2} (D_{2,1}D_{3,1})^{-2}(D_{3,1}-D_{2,1})^2 \sin^2 \theta \cos^2 \theta 
\right] \mx{,}
\end{aligned}
\eeq
where $\tilde{\mathfrak{W}}$ stands for $\mathfrak{W}$ in the new set of variables.  It is evident that for sufficiently small $\ep,\mu$ and all $\alpha>0$, the latter is strictly negative. Our goal is now to find a bound for the thresholds $\ep_0,\mu_0$ in such a way this property persists for all $\ep \leq \ep_0$ and $\mu \leq \mu_0$. \\
Let us firstly deal with $O(\ep^2)$ coefficients. For this purpose, it is immediate to check that 
\[
1+2 \mu x(\sin^2 \theta + \gamma^{-1} \cos^2 \theta)+ (\mu x)^2(\sin^2 \theta + \gamma^{-2} \cos^2 \theta) \geq 1-2 \mu/\gamma \mx{.}
\]
Hence the latter will be, say, greater than $1/2$, for all $\mu\leq \mu_0:=\gamma/4$.\\
As for the terms of $O(\ep^3)$, by recalling (\ref{eq:propertiesd}), one has
\[
 (D_{2,1}^{-1}-D_{3,1}^{-1})\cos \theta \sin^2 \theta \leq D_{2,1}^{-1} \mx{.}
\] 
Finally, $O(\ep^4)$ contributions are clearly strictly positive. In conclusion one has
\[
\dot{\tilde{\mathfrak{W}}} < \alpha \ep^2 \left[2 \Omega \ep /  (D_{2,1} \gamma) -1/2 \right]\mx{.}
\]    
The latter implies that the derivative along the solutions $\tilde{\mathfrak{W}}$ will be strictly negative by choosing
\[
\ep_0:=D_{2,1} \gamma/(4 \Omega) \mx{.}
\]
This implies that, in the original set of variables $\ue{U}$, the cylinder ellipse - shaped region on which $\dot{\mathfrak{W}}$ is strictly negative, is given by
\[
\ml{C}_{el}:=\{\ue{U}\in \RR^3:|U_1|\leq \gamma/4, \quad  U_2^2+(U_3 D_{3,1}/(D_{2,1} \gamma))^2=(4 \Omega)^{-2} \} \mx{.}
\] 
In order to obtain the required estimate for the basin of attraction, we need to find the largest value $\mathfrak{W}^*$ for which the level curve $\mathfrak{W}(U_2,U_3)=\mathfrak{W}^*$ is entirely contained in the projection of $\ml{C}_{el}$ on the $(U_2,U_3)$ plane. This is a simple problem of constrained optimisation in which $\mathfrak{W}$ has to be minimised on the ellipse $U_2^2+(U_3 D_{3,1}/(D_{2,1} \gamma))^2=(4 \Omega)^{-2}$. It is easy to check that the required minimum is attained at 
\[
(U_2^*,U_3^*):=(0,\pm D_{2,1}\gamma/(4 D_{3,1} \Omega)) \mx{,}
\]   
implying $\mathfrak{W}^*:=\mathfrak{W}(U_2^*,U_3^*) \equiv (D_{2,1} \gamma)^2/(32 D_{3,1} \Omega^2)$. It is now necessary to check that $|U_1(t)|\leq \gamma/4$ for all $t$, ensuring that $\ue{U}(t) \in \ml{C}_{el}$ and then $\dot{\mathfrak{W}} <0$ all $t$. \\
For this purpose, let us observe that the level $\tilde{\Psi}(\ue{U})=\ml{E}$ is equivalent to
\[
(U_1+\gamma)^2+(U_2-\Omega)^2+U_3^2=1+\delta_a
\]
The latter defines a $\delta_a$-family of surfaces in which $U_1$ can be written as a graph (surface) over the other two variables in a suitable neighbourhood of the origin, as follows   
\beq{eq:surfaceuone}
U_1=U_1(U_2,U_3;\delta_a) \equiv -\gamma+\sqrt{1+\delta_a-(U_2-\Omega)^2-U_3^2} \mx{,}
\eeq
and such that $U_1(0,0;0)=0$. The sign has been chosen according to the fact that we are describing the portion of sphere located near $\ue{s}^+$. The task consists in showing that, under the required condition for $\delta_a$, any point $(U_2,U_3) \in \ml{B}_{\ml{E}}$, yields $|U_1| \leq \gamma/4$. Given the structure of (\ref{eq:surfaceuone}), it will be sufficient to show $U_1\geq-\gamma/4$, being condition $U_1\leq \gamma/4$ a consequence. From (\ref{eq:surfaceuone}), one has $U_1 \geq -\gamma + \sqrt{-|\delta_a|+\gamma^2-A_a}$, where $A_a:=2 \Omega |U_2|+U_2^2+U_3^2$. It is immediate to verify that, by choosing $\delta_a$ as in (\ref{eq:conddeltaa}), the desired condition is obtained if 
\beq{eq:conditionaa}
A_a \leq (3/8) \gamma^2 \mx{.} 
\eeq
On the other hand, in $\ml{B}_{\ml{E}}$ one has 
\[
|U_2| \leq [\gamma/(4 \Omega)] q_a, \qquad U_2^2+U_3^2 \leq [\gamma/(4 \Omega)]^2 q_a^2 \mx{,}  
\] 
where $q_a:=\sqrt{D_{2,1}/D_{3,1}}$. This implies that condition (\ref{eq:conditionaa}) holds if
\[
4 q_a + 2^{-1} \Omega^{-2} q_a^2  \gamma \leq \gamma \mx{.}
\] 
A sufficient condition for the latter to hold is that, for instance, $4 q_a \leq \gamma/4$ and $2^{-1} \Omega^{-2} q_a^2 \leq 3/4 $. They are equivalent to the first and the second of  (\ref{eq:choiceomegak}), respectively. Equation (\ref{eq:asymptotic}) follows directly from (\ref{eq:surfaceuone}) whilst bound (\ref{eq:absu}) is immediate from (\ref{eq:asymptotic}) and a Taylor estimate.\\
As for the last statement it is sufficient to observe that the ellipse resulting from the projection of $\ml{B}_{\ml{E}}$ clearly contains the disk whose radius is not larger than the minor semi-axis of the former, which is given by the r.h.s. of (\ref{eq:semiminor}).
\endproof
\subsubsection*{Proof of Thm. \ref{thm:main}}
\proof Let us consider the class of expulsion targets and currents such that $\beta_e,K^2=O(\lambda)$. 
In order to use Lemmata \ref{lem:expulsion}, \ref{lem:transfer} and \ref{lem:attraction} we are going to make the assumptions required by these results. First of all, given $D_{2,1}=O(\lambda)$ and $D_{3,1}=O(1)$, let us choose $\Omega$ within the limitation prescribed by the first of (\ref{eq:choiceomegak}), this determines a range for $h_2=-D_{2,1}\Omega$. It is easy to realise that $\Omega=O(\sqrt{\lambda})$. As a consequence, $\gamma\sim 1-O(\lambda)$, hence we can assume without loss of generality, that the last of (\ref{eq:choiceomegak}) is satisfied as well. There is a certain amount of freedom in the choice of $D_{3,2}$, provided that (\ref{eq:propertiesd}) holds. 
We anticipate that the value of $\delta_a$ will be chosen, but only at the very end of the proof, according to (\ref{eq:conddeltaa}), hence we shall proceed by supposing that such a condition is satisfied.\\
Once the possibility to use Lem. \ref{lem:attraction} has been guaranteed, let us examine the setting of the remaining Lemmata. As for Lem. \ref{lem:expulsion}, the values of $r^*,\ml{M}_{1,2}$ are uniquely determined, this yields $\lambda_e$. The same holds for the constants of $
\ml{M}_{tr},\ml{K}_{\ue{w}}$ and the threshold $\lambda_{tr}$ of Lem. \ref{lem:transfer}. In order to proceed with the proof we shall assume $\lambda_0:=\min\{\lambda_e,\lambda_{tr}\}$, although practical applications of Lem. \ref{lem:transfer} may require a further restriction of such a threshold, as already discussed in Rem. \ref{rem:threshold}. \\
Let us now use the mentioned Lemmata to construct a suitable neighbourhood of the equilibrium $\ue{s}^-$ which is transported in the vicinity of $\ue{s}^+$ and then attracted by the point close to it on the corresponding energy level. The proof will be complete once (\ref{eq:transferstat}) is validated and a precise value for  $r^-$ is computed. \\
For this purpose let us start by observing that, from (\ref{eq:apprxitwo}) and (\ref{eq:translation}), at the end of the expulsion stage started in $\ue{u}_c(0)\equiv \ue{s}^-$, we get 
\[
\ue{u}_c(T_e)=(-\gamma +\lambda \kappa \bar{a},-\Omega - \lambda \kappa \bar{b},-K ) \mx{,}
\]
with $\kappa:=(2 \beta_e)^{-1}(-1+\sqrt{1-2 K^2})\sim K^2/(2 \beta_e)=O(1)$ by assumption and $T_e$ given by (\ref{eq:texp}). Note that $T_e \sim K/\beta_e = O(\lambda^{-1/2})$. Furthermore, $\Phi(\ue{u}_c(T_e))=1+O(\lambda)$. \\
Let us now set $\ue{w}_r(0):=\ue{u}_c(T_e)$ in such a way to start the transfer stage. This determines the value of $T_r$ via (\ref{eq:ttr}). By construction we have $w_1^{[0]}(T_{tr})=-\ml{A}_m-K^2=\gamma+O(\lambda)=:\hat{\gamma}$ and $w_3^{[0]}(T_{tr})=-K$. In order to determine $w_2^{[0]}(T_{tr})$, one can use the conservation law (\ref{eq:firstint}), which reads as 
\[
\sigma^{-1}(\lambda \kappa \bar{a}-\gamma)^2+\sigma(\Omega+\lambda \kappa \bar{b})^2=\sigma^{-1}\hat{\gamma}^2+\sigma (w_2^{[0]}(T_r))^2 \mx{.} 
\]
The latter gives, recalling (\ref{eq:abarbbar}) and expanding in $\lambda$,  
\beq{eq:wrtwo}
w_2^{[0]}(T_r)=-\Omega+O(\lambda) \mx{.}
\eeq
The sign in the r.h.s. of the latter has been chosen by using the fact that before and on the first occurrence of 
$w_1^{[0]}=-\ml{A}_m-K^2$, the trajectory lies in the $w_2^{[0]}$-negative half-plane. Hence, by recalling (\ref{eq:integrationv}), we have
\beq{eq:wte}
\ue{w}(T_e)=(\gamma,-\Omega,-K)+O(\lambda) \mx{.}
\eeq
Once more, the usual conservation law is easily checked as $\Phi(\ue{w}(T_e))=1+O(\lambda)$.\\
Now we need to make sure that (\ref{eq:wte}) is actually contained in the basin of attraction $\ml{B}_{\ml{E}}$ in such a the attractivity property stated in Lem. \ref{lem:attraction} can be used. For this purpose, we recall the very last statement of Lem. \ref{lem:attraction} and set $r_{sm} \geq (5/4)K$ (i.e. slightly bigger than $K$ itself). Hence, (\ref{eq:semiminor}) yields 
\[
\Omega^2 \leq (5 K D_{3,1})^2 (\gamma D_{2,1})^2 \mx{,}
\]
yielding, by recalling that $\gamma^2=1-\Omega^2$, $\Omega^2 \leq [1+25 (D_{3,1}/D_{2,1}^2)K^2]^{-1}$. \\The latter, if compared with the first of (\ref{eq:choiceomegak}) (which provides a lower bound for $\Omega^2$), makes sense provided that $50 K^2 D_{3,1}\leq 3 \gamma^2 D_{2,1} $. One could choose, for instance, 
\[
K=\bar{K}:=(\gamma/4) \sqrt{D_{2,1}/D_{3,1}} \mx{,}
\]
then computing $\Omega^2$ as the mid-point given by the two conditions
\[
\Omega^2:=2^{-1}[(2/3)(D_{2,1}/D_{3,1})+\gamma^2/(25 \bar{K}^2)(D_{2,1}/D_{3,1})^2]=(52/75)(D_{2,1}/D_{3,1}) \mx{,}
\]
this determines $h_2$. \\
Hence, any suitable manifold $\ml{E}=1+\delta_a$ contained in the cylinder $\{\sqrt{U_2^2+U_3^2}\leq (5/4)K,\, \, |U_1|\leq \gamma/4 \}$ is contained in $\ml{B}_{\ml{E}}$. This implies that any 
$\mathfrak{B}_{K/8}(\ue{w}(T_e))$ is contained in $\ml{B}_{\ml{E}}$ as well. Hence, we can set $\rho_{tr}^+:=K/(8 \lambda)$ and $\lambda \rho_{tr}^-:=[32(1+\Theta)]^{-1}\Theta K$ by (\ref{eq:choicerhotr}), then 
use Lem. \ref{lem:transfer}. The latter, allows us to conclude that the set $\mathfrak{B}_{\lambda \rho_{tr}^-}(\ue{w}(T_e))$ will evolve inside $\ml{B}_{\ml{E}}$ under the action of the controlled phase flow $\Phi_{tr}^t$.\\
As $\omega=O(\lambda^{-1/2})$ because of the choice for $K$, we have $\Theta=O(\sqrt{\lambda})$, hence $\lambda \rho_{tr}^-=O(\lambda)$. It is now sufficient to recall Lem. \ref{lem:expulsion} and set $\rho_e:=\rho_{tr}^-$ to determine the value of $r^-$, which is easily found via (\ref{eq:lamz}) as
\[
r^-=[128 \ml{M}_e (1+\Theta)]^{-1}\Theta K \mx{.}
\]
Similarly, we notice that $r^-=O(\lambda)$. \\ 
As (\ref{eq:absu}) holds for the family of points starting in $\ml{B}_{\ml{E}}$ and parameterised by $\delta_a$ in the admissible range (\ref{eq:conddeltaa}), this property will be true, \emph{a fortiori}, for the points of the set 
\[
\ml{U}_a:=\Phi_{tr}^{T_{tr}}\left(\Phi_e^{T_e}(\mathfrak{B}_{r^-}(\ue{s}^-)) \right) \subset \ml{B}_{\ml{E}} \mx{.}
\]
Hence, by defining
\[
\ue{u}_{\infty}:=\lim_{t \rw + \infty} \Phi_a^t \left( \bigcup_{\ml{E}=1+\delta_a } \ml{B}_{\ml{E}} \right) \mx{,}
\]
its distance from $\ue{s}^+$ is readily bounded by $\delta_a$ via (\ref{eq:absu}) as follows
\[
|\ue{u}_{\infty}(\delta_a)-\ue{s}^+|=|\ue{U}_{\infty}(\delta_a)| \leq \delta_a/(2 \gamma) \mx{.}
\]
On the other hand, by the conservation of $\Psi(\ue{u})$, an upper bound for $\delta_a$ is found by evaluating the maximum of $|\Psi(\ue{u})-1|$ for all $\ue{u}:=\ue{s}^- + \ue{u}' \in \mathfrak{B}_{r^-}(\ue{s}^-)$ i.e. for all $\ue{u}' \in \mathfrak{B}_{r^-}(\ue{0})$. This is straightforward, as $\Psi(\ue{s}^- + \ue{u}')\leq 1 + r^-$ by the triangle inequality, hence $|\delta_a| \leq r^- $. This proves (\ref{eq:transferstat}) with $f:=(2 \gamma)^{-1}$. The proof of Thm. \ref{thm:main} is now complete. 
\endproof
\begin{figure*}[t!]\begin{center}
		\begin{minipage}[c][1\width]{0.48\textwidth}
			\hspace{0pt} 
			{\begin{overpic}[width=\textwidth]{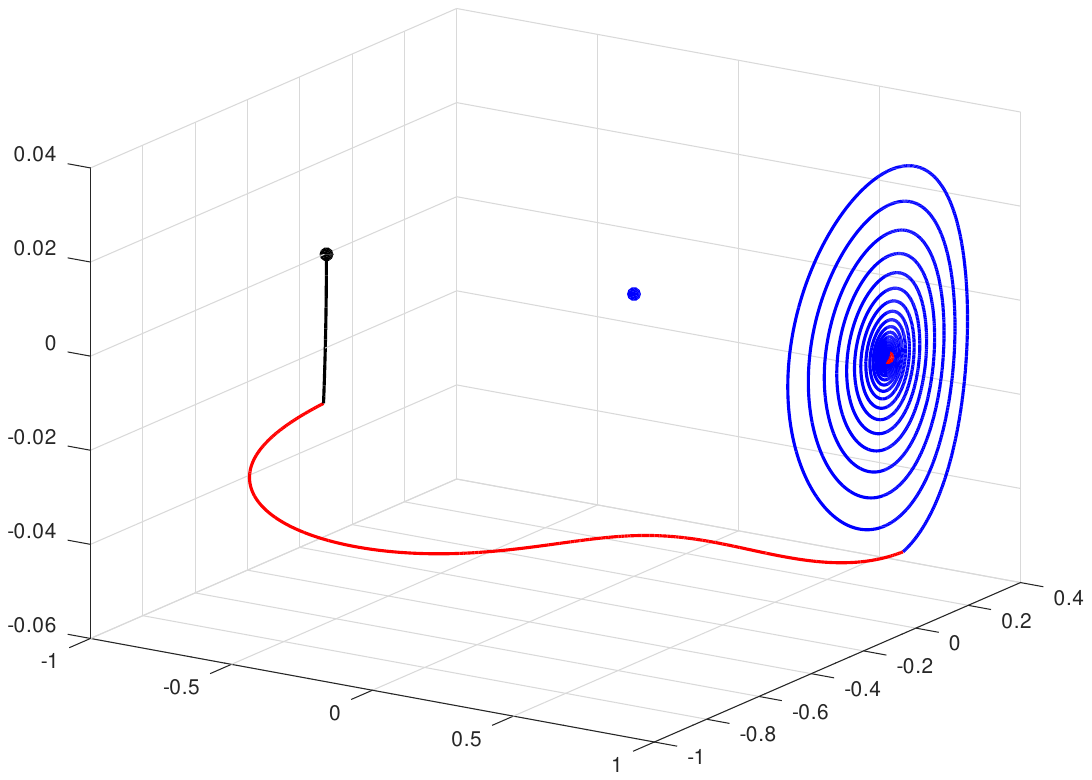}
					\put(0,78){\footnotesize (a)}
   		 		\put(25,60){\footnotesize $\ue{s}^{-}$}
         	        \put(44,58){\footnotesize $\ue{0}$}
                  	\put(2,58){\footnotesize $u_3$}
     	            \put(17,28){\footnotesize $u_1$}
					\put(57,28){\footnotesize $u_2$}
             \end{overpic}}
		\end{minipage}
    		\begin{minipage}[c][1\width]{0.48\textwidth}
			\hspace{-5pt}
			{\begin{overpic}[width=1\textwidth]{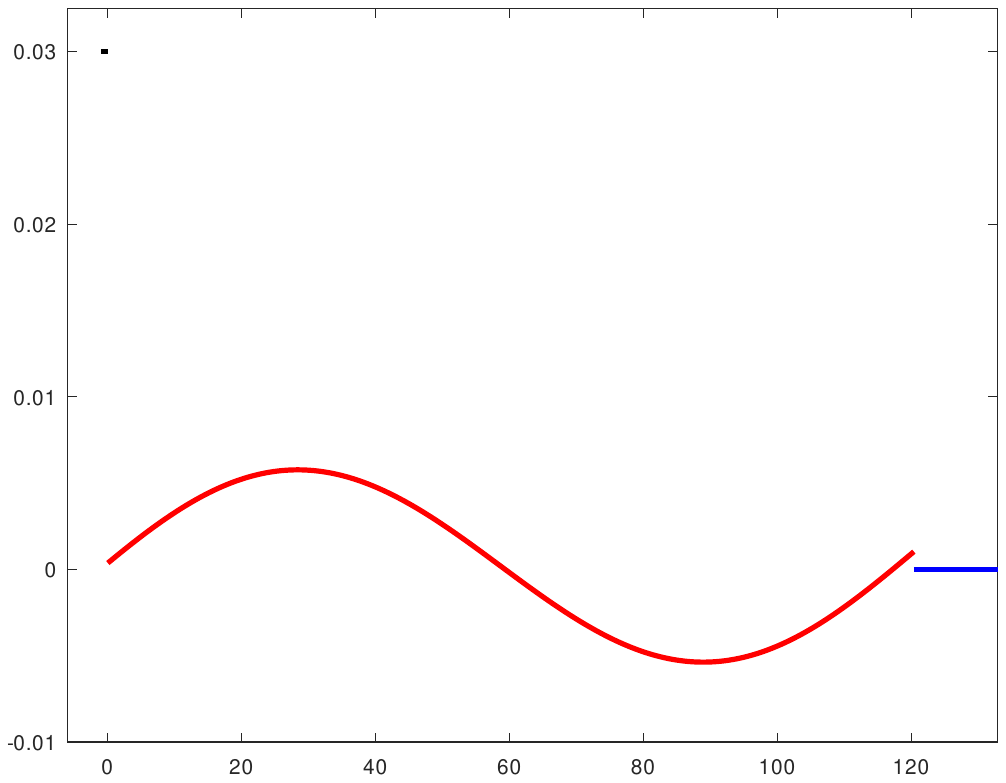}
					\put(0,78){\footnotesize (b)}
					\put(16,69.75){\footnotesize $\longleftarrow$}
        	     	\put(5,65){\footnotesize $\beta(t)$}
					\put(52,23){\footnotesize $t$}
 			\end{overpic}}
		\end{minipage}	\\
  		\vspace{-30pt}
	\end{center}
	\vspace{-20pt}
	\caption{Example of a full switching process. Panel (a) shows the (piece-wise) trajectory in the space $\ue{u}$, as a sequence of an ``expulsion'' (black), a ``transfer'' (red) and finally the ``attractive'' stage (blue). The initial condition has been chosen as $\ue{u}(-T_e)=\ue{s}^-+\lambda (-0.1,0.05,0)$, with $T_e=1.0262$. The first $\sim 133$ time units of the constructed controlling current are reported in panel (b). The small black segment for $t\in[-T_e,0]$ (emphasised by the arrow), represents the current injected during the expulsion stage. We have set $\lambda=0.002$. The full set of parameters is reported in Appendix C.}
	\label{fig:full}
	\vspace{0pt}
\end{figure*}

\section{Validation tests and conclusions}
\begin{figure*}[t!]\begin{center}
		\begin{minipage}[c][1\width]{0.48\textwidth}
			\hspace{0pt} 
			{\begin{overpic}[width=\textwidth]{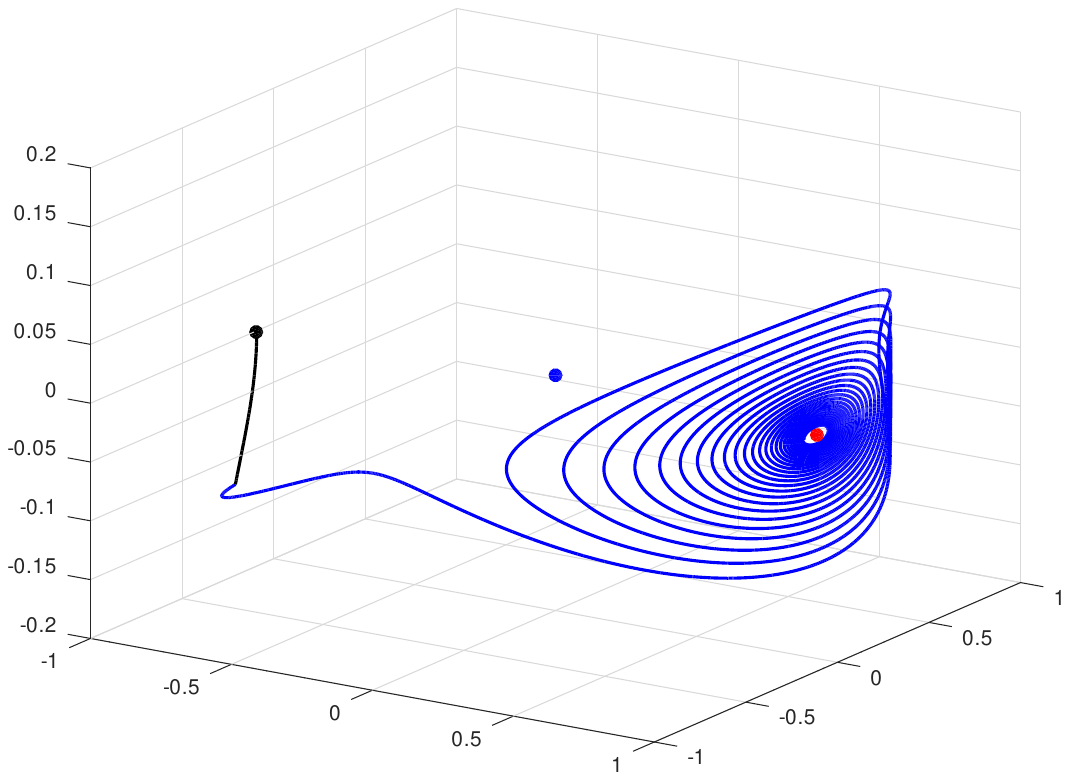}
					\put(0,78){\footnotesize (a)}
   					\put(2,58){\footnotesize $u_3$}
     	            \put(17,28){\footnotesize $u_1$}
					\put(57,28){\footnotesize $u_2$}
             \end{overpic}}
		\end{minipage}
    		\begin{minipage}[c][1\width]{0.48\textwidth}
			\hspace{-5pt}
			{\begin{overpic}[width=1\textwidth]{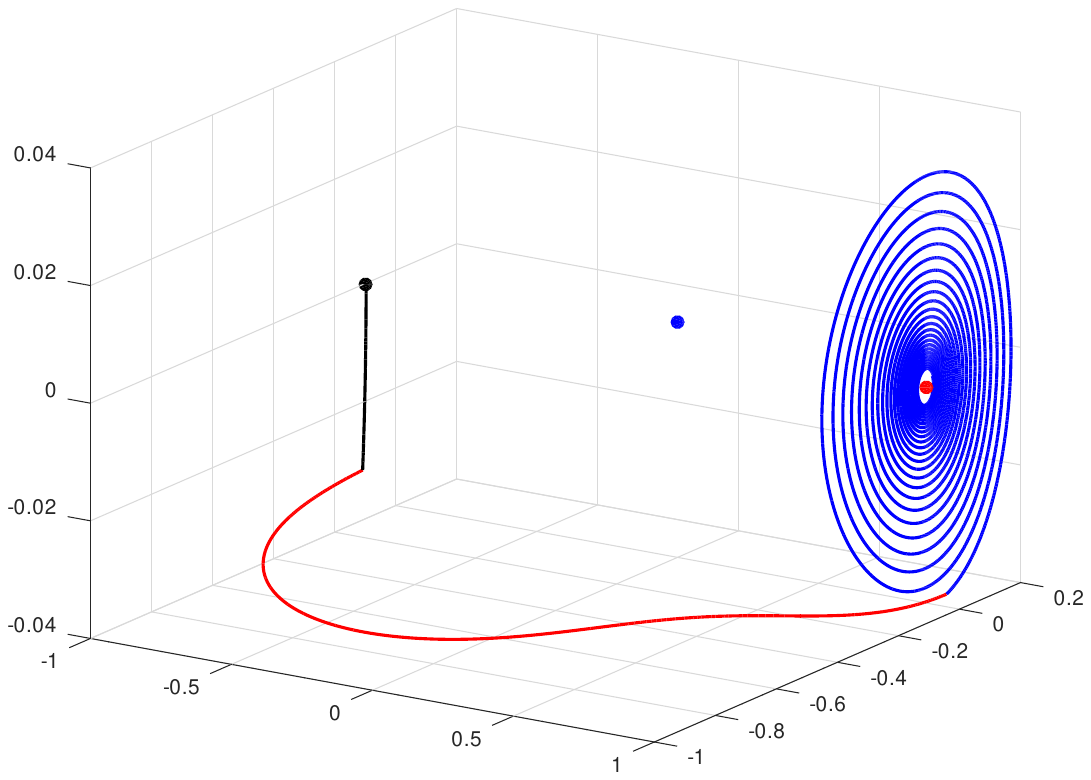}
					\put(0,78){\footnotesize (b)}
 					\put(2,58){\footnotesize $u_3$}
     	            \put(17,28){\footnotesize $u_1$}
					\put(57,28){\footnotesize $u_2$}
  			\end{overpic}}
		\end{minipage}	\\
  		\vspace{-30pt}
  	\begin{minipage}[c][1\width]{0.48\textwidth}
			\hspace{5pt} 
			{\begin{overpic}[width=\textwidth]{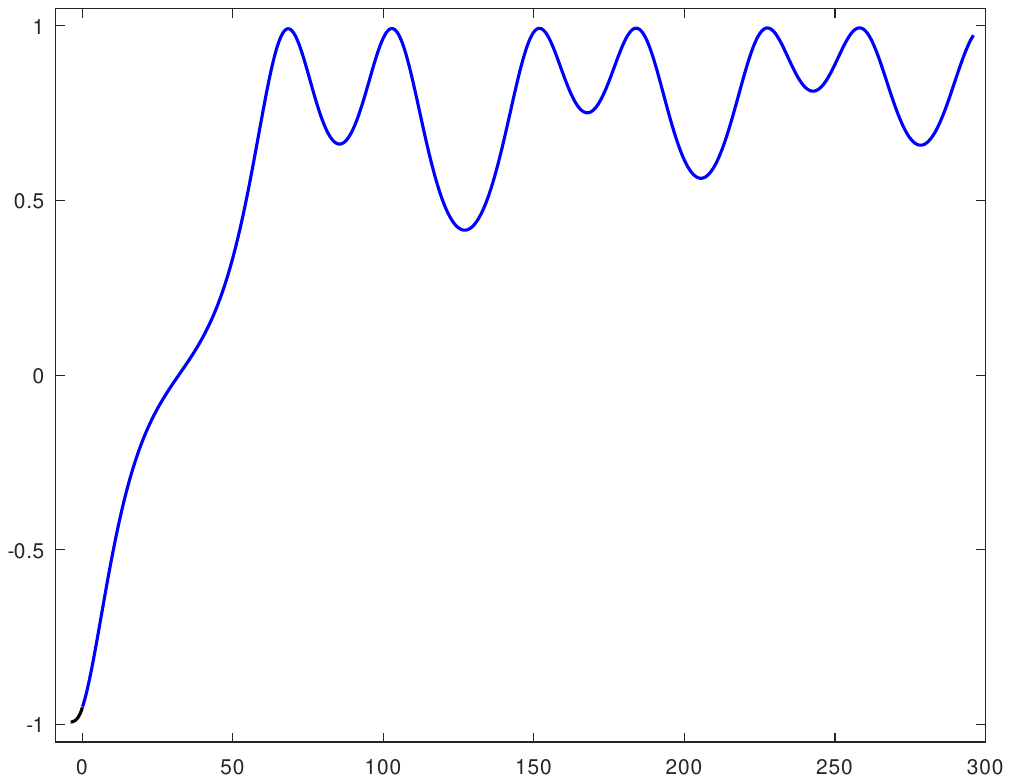}
					\put(0,78){\footnotesize (c)}
 					\put(4,66){\footnotesize $u_1$}
 					\put(54,24){\footnotesize $t$}
			\end{overpic}}
		\end{minipage}	
    \begin{minipage}[c][1\width]{0.48\textwidth}
			\hspace{0pt}
			{\begin{overpic}[width=\textwidth]{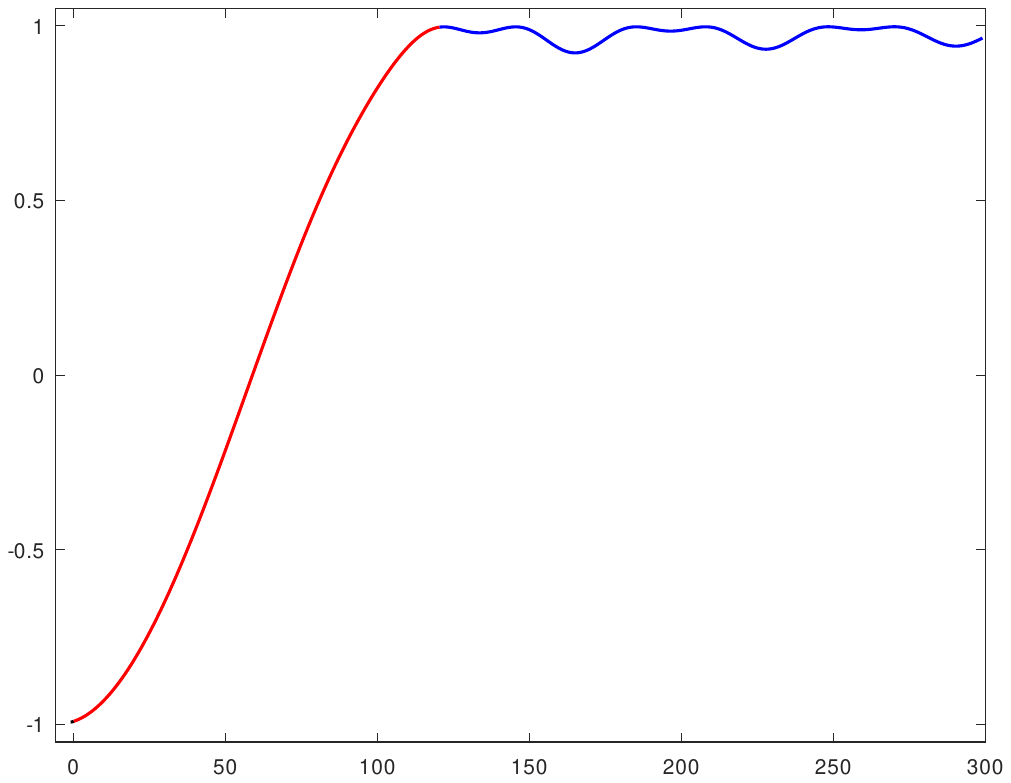}
					\put(0,78){\footnotesize (d)}
					\put(4,66){\footnotesize $u_1$}
					\put(54,24){\footnotesize $t$}
			\end{overpic}}
		\end{minipage}\\
	\end{center}
	\vspace{-30pt}
	\caption{Comparison between the ballistic switching (a) and the one obtained by using CQL solutions (b). The initial condition has been set as $\ue{u}(-T_e)=\ue{s}^-+(\lambda,0,0)$ in both cases. Below, ((c) and (d), respectively), the corresponding behaviour of the first coordinate $u_1(t)$. A considerably small residual oscillation of the mentioned coordinate can be observed in the CQL case. The remaining parameters used are specified in Appendix C.}
	\label{fig:comparison}
	\vspace{0pt}
\end{figure*}

\begin{figure*}[t!]\begin{center}
		\begin{minipage}[c][1\width]{0.48\textwidth}
			\hspace{0pt} 
			{\begin{overpic}[width=\textwidth]{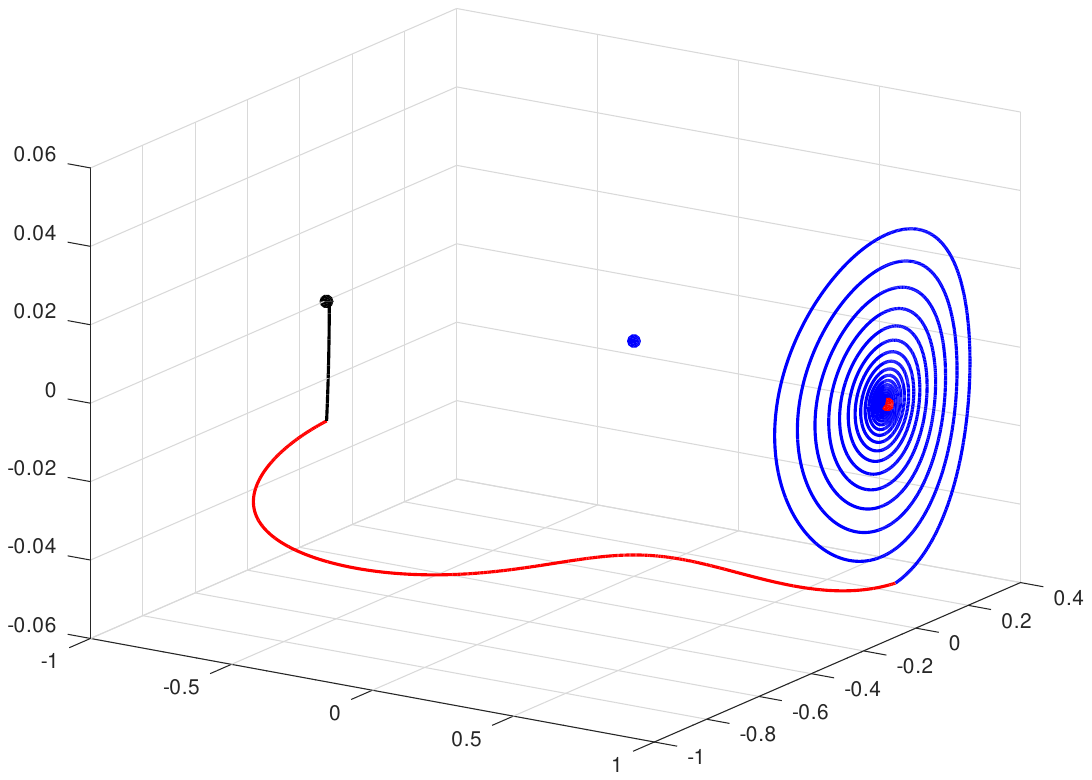}
					\put(0,78){\footnotesize (a)}
 					\put(2,58){\footnotesize $u_3$}
     	            \put(17,28){\footnotesize $u_1$}
					\put(57,28){\footnotesize $u_2$}
             \end{overpic}}
		\end{minipage}
    		\begin{minipage}[c][1\width]{0.48\textwidth}
			\hspace{-5pt}
			{\begin{overpic}[width=1\textwidth]{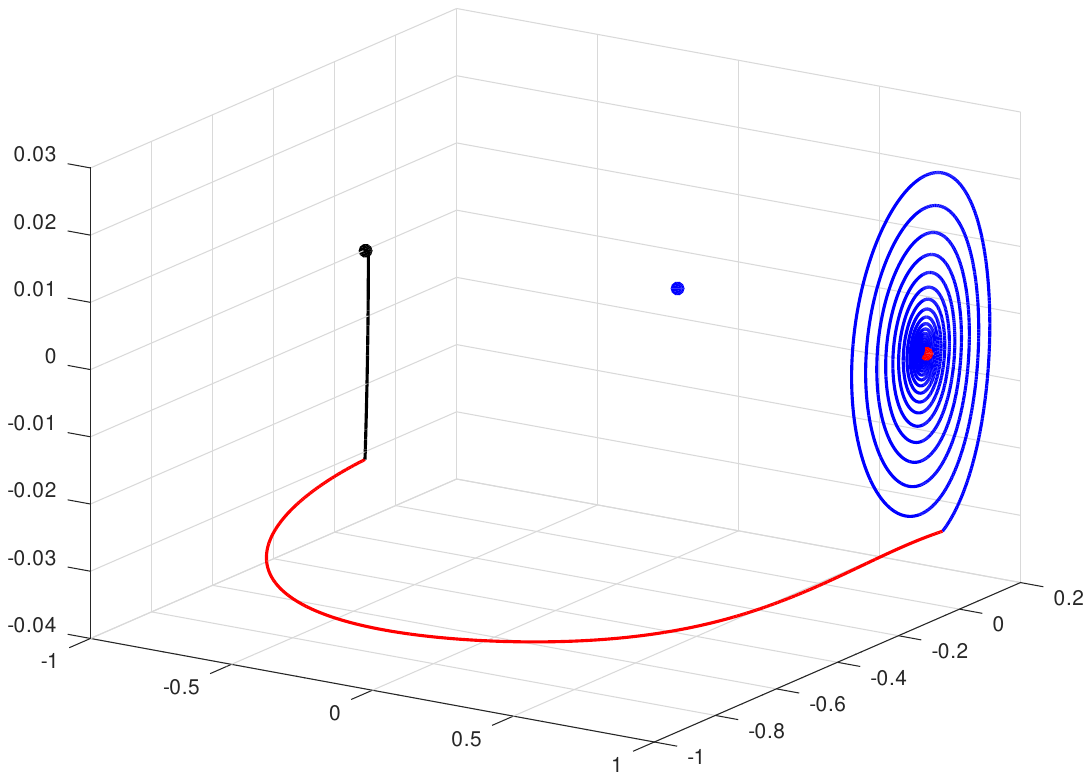}
					\put(0,78){\footnotesize (b)}
 					\put(2,58){\footnotesize $u_3$}
     	            \put(17,28){\footnotesize $u_1$}
					\put(57,28){\footnotesize $u_2$}
			\end{overpic}}
		\end{minipage}	\\
  		\vspace{-30pt}
  	\begin{minipage}[c][1\width]{0.48\textwidth}
			\hspace{5pt} 
			{\begin{overpic}[width=\textwidth]{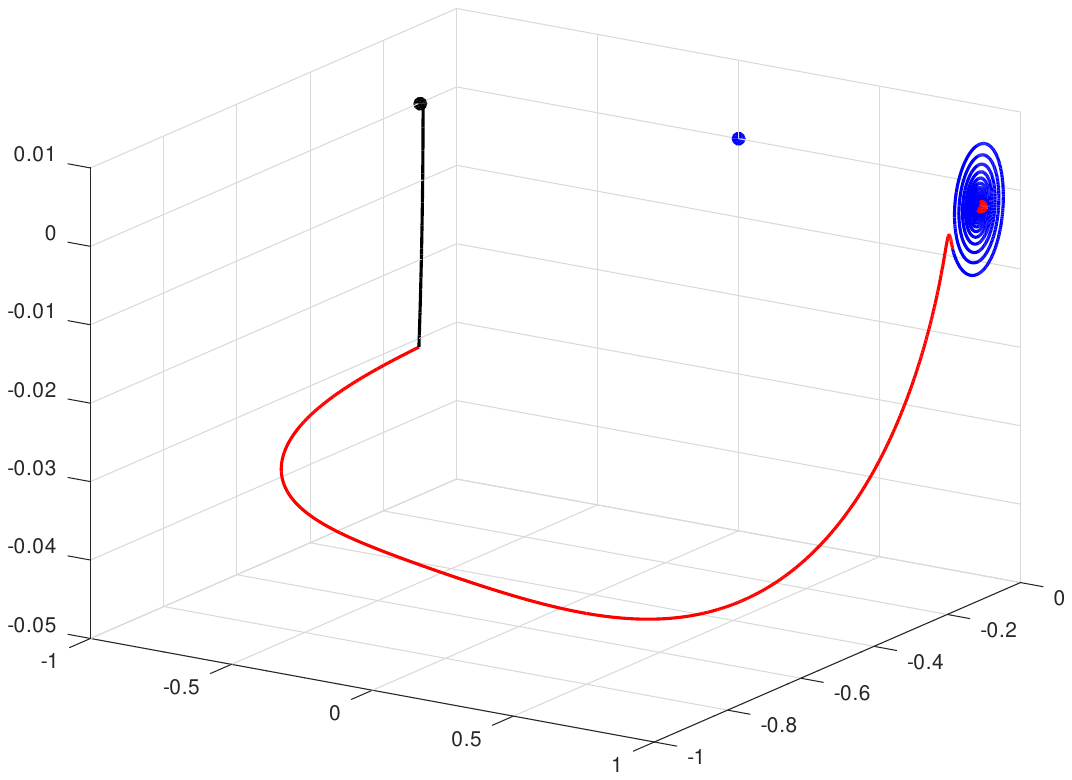}
					\put(0,78){\footnotesize (c)}
 					\put(2,58){\footnotesize $u_3$}
     	            \put(17,28){\footnotesize $u_1$}
					\put(57,28){\footnotesize $u_2$}
			\end{overpic}}
		\end{minipage}	
    \begin{minipage}[c][1\width]{0.48\textwidth}
			\hspace{0pt}
			{\begin{overpic}[width=\textwidth]{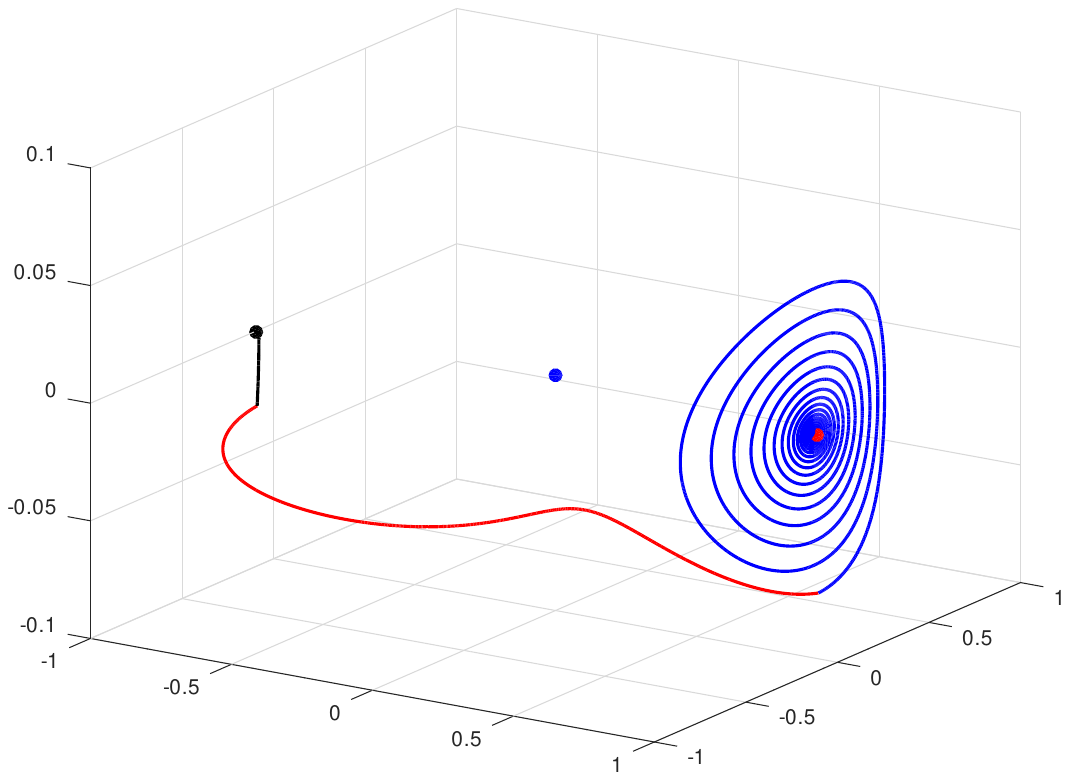}
					\put(0,78){\footnotesize (d)}
 					\put(2,58){\footnotesize $u_3$}
     	            \put(17,28){\footnotesize $u_1$}
					\put(57,28){\footnotesize $u_2$}
			\end{overpic}}
		\end{minipage}\\
	\end{center}
	\vspace{-30pt}
	\caption{Switching trajectories via the CQL approach with a simulated error in the control: panels (a) and (b) show the case of a time contraction/dilatation of the computed ``espulsion'' time,  more precisely $T_e \leftarrow 0.98 T_e$ and $T_e \leftarrow 1.02 T_e$, respectively. Panels (c) and (d) report an analogous test during the ``transfer'' stage instead. Hence, the time variable in the control $\tilde{\beta}_{tr}(t)$ is replaced with $t \leftarrow  0.98 t$ and $t \leftarrow 1.02 t$, respectively. Taking into account the ``considerable'' duration of the transfer process (it has an order of $10^2$ time units) the latter should be regarded as a pretty ``hard'' test and evident trajectory modifications such as those shown in (c) and (d), are reasonably expected. Full set of parameters in Appendix C.}
	\label{fig:stress}
	\vspace{0pt}
\end{figure*}

The aim of this section is to provide examples of full switching processes realised according to Thm. \ref{thm:main}. To summarise, after an initial ``expulsion'' stage, the evolving point is ``transferred'' via a CQL solution in the estimated basin of attraction of the target equilibrium so that it can be ``attracted'' by it. An example of a full switching is reported in fig. \ref{fig:full}.\\
As anticipated in sec. \ref{sec:one}, the approach proposed here is conceptually different from the well established ballistic switching procedure, see, e.g. \cite{10.1063/1.4914122}. From the qualitative viewpoint, these stages could be compared to our  ``expulsion'' and ``attraction'' respectively, see fig. \ref{fig:comparison}, panel (a). However, as it is evident from the latter, in order to successfully complete the switching, the operator needs to trust a basin of attraction which is, for instance, way larger than the one rigorously computed in this work. Clearly, the one stated in Lem. \ref{lem:attraction} is nothing but a sufficient condition but it suggests that, if not addressed with further (and highly specialised) tools and ad-hoc arguments, the attractivity of larger sets retains no more than a probabilistic validity. On the contrary, the proposed CQL strategy has a genuinely deterministic character, provided that $\lambda$ is sufficiently small (according to the bounds described in the proof) and the initial set of initial conditions is chosen in $\mathfrak{B}_{r^-}(\ue{s}^-)$. \\
In addition, as anticipated in sec. \ref{sec:one}, our method only requires a minimal amount of initial injected current during the ``expulsion'' stage, as the solution is subsequently ``guided'' towards $\ml{B}_{\ml{E}}$ via the constructed control, along a CQL trajectory. As a result, this remarkably reduces the ``energy'' carried by $u_3$ and then the amplitude of the oscillations around $\ue{s}^+$ when the current is switched off. See fig. \ref{fig:comparison} for a depiction of this phenomenon in a comparison with the ballistic approach. \\
As a further check in relation with the robustness of this method and its deterministic character, the outcomes of a stress test are proposed in fig. \ref{fig:stress}, in which some errors are simulated with respect to the proposed control. More precisely, either a time dilatation-contraction coefficient for the expulsion time $T_e \leftarrow j T_e$ with $j=1\pm 0.02$ or for the transfer control, i.e. $\title{\beta}_{tr}(t) \leftarrow \tilde{\beta}_{ tr}(j t)$ are considered. \\
\medskip \\
In conclusion, the present work proposes an analytical formulation apt to realise a fully deterministic switching mechanism. The argument relies on the concept of CQL solution, a highly non-local object at the heart of the procedure, constructed via a perturbative approach which exploits intrinsic technologically relevant features of the system. The mentioned CQL solutions are realised by the explicit determination of a time dependent control, determined by using well established tools of Hamiltonian perturbation theory. Interestingly, the latter possess the potential to increase the range of validity of the method to a even larger class of systems, should a higher order analysis be considered. The numerical experiments proposed offer either a validation or a visual interpretation of the main statement of Thm. \ref{thm:main} and the related Lemmata.

\section*{Appendix A}
Explicit expressions of $\ue{\ml{F}}(\ue{u};\ue{v})$ and $\ue{\ml{R}}(\ue{u};\ue{v})$:
\[
\begin{aligned}
\ml{F}_1&=-u_1 u_3 (-v_2^2 \tilde{\alpha} D_{3,2} K-v_1^2 \tilde{\alpha} D_{3,1} K+v_1 \tilde{h}_2-v_1 v_2 \tilde{D}_{2,1}) \rho+u_3 \tilde{h}_2+u_1 u_3^2 \tilde{\alpha} D_{3,1}\\
\ml{F}_2&=u_2 u_3^2 \tilde{\alpha} D_{3,2}-u_2 u_3 (-v_2^2 \tilde{\alpha} D_{3,2} K-v_1^2 \tilde{\alpha} D_{3,1} K+v_1 \tilde{h}_2-v_1 v_2 \tilde{D}_{2,1}) \rho\\
\ml{F}_3&=(u_2^2+u_1^2) (-v_2^2 \tilde{\alpha} D_{3,2} K-v_1^2 \tilde{\alpha} D_{3,1} K+v_1 \tilde{h}_2-v_1 v_2 \tilde{D}_{2,1}) \rho-u_1 \tilde{h}_2\\
&+\tilde{\alpha} (-u_2^2 u_3 D_{3,2}-u_1^2 u_3 D_{3,1})+u_1 u_2 \tilde{D}_{2,1}\\
\ml{R}_1&=u_1 v_2 u_3 \tilde{\alpha} \tilde{h}_2 K \rho+\tilde{\alpha} (u_1 u_2^2 \tilde{D}_{2,1}-u_1 u_2 \tilde{h}_2)\\
\ml{R}_2&=u_2 v_2 u_3 \tilde{\alpha} \tilde{h}_2 K \rho+\tilde{\alpha} (u_3^2 \tilde{h}_2+u_1^2 \tilde{h}_2-u_1^2 u_2 \tilde{D}_{2,1})\\
\ml{R}_3&=-(u_2^2+u_1^2) v_2 \tilde{\alpha} \tilde{h}_2 K \rho-u_2 u_3 \tilde{\alpha} \tilde{h}_2
\end{aligned}
\]

\section*{Appendix B}
Explicit expression of the entries of $\{\uue{D} \ue{V}(\ue{\xi},\tilde{\beta})\}$. Let us denote them with $\lambda b_{i,j}$. We have 
\[
\begin{aligned}
b_{1,1}&=\xi_3\tilde{\beta}\lambda + \tilde{\alpha} (\tilde{D}_{2,1}\xi_2\Omega+D_{3,1}\xi_3^2)\lambda^2 + \tilde{\alpha} \tilde{D}_{2,1}\xi_2^2\lambda^3 \\
b_{1,2}&=D_{3,2} \xi_3 -\tilde{\alpha}\tilde{D}_{2,1}\Omega \thd \lambda  
+\tilde{\alpha} \tilde{D}_{2,1} (\xi_1 \Omega-2 \xi_2 \thd)\lambda^2+
2\tilde{\alpha}\tilde{D}_{2,1} \xi_1 \xi_2 \lambda^3 \\
b_{1,3}&=(-\tilde{\beta} \thd +\tilde{D}_{2,1}\Omega+D_{3,2}\xi_2)+(\xi_1\tilde{\beta}-2\tilde{\alpha} D_{3,1} \xi_3 \thd )\lambda+2\tilde{\alpha}D_{3,1}\xi_1\xi_3\lambda^2 \\
b_{2,1}&=-D_{3,1} \xi_3 +2\tilde{\alpha}\tilde{D}_{2,1}\xi_2 \thd \lambda^2 -2\tilde{\alpha}\tilde{D}_{2,1}\xi_1\xi_2\lambda^3 \\
b_{2,2}&=(\xi_3\tilde{\beta}-\tilde{\alpha}\tilde{D}_{2,1}\gamma^2)\lambda + \tilde{\alpha} (2 \tilde{D}_{2,1}\xi_1 \thd +D_{3,2}\xi_3^2)\lambda^2 -\tilde{\alpha}\tilde{D}_{2,1}\xi_1^2\lambda^3 \\
b_{2,3}&=\tilde{\beta}\Omega-D_{3,1}\xi_1 +(2\tilde{\alpha} D_{3,2} \xi_3\Omega+\xi_2\tilde{\beta})\lambda+ 2\tilde{\alpha} \xi_3 (\tilde{D}_{2,1} \Omega+D_{3,2} \xi_2)\lambda^2 \\
b_{3,1}&=2\tilde{\beta} \thd (2\tilde{\alpha}D_{3,1} \xi_3 \thd-2\xi_1\tilde{\beta}+\tilde{D}_{2,1}\xi_2)\lambda -2\tilde{\alpha}D_{3,1} \xi_1\xi_3\lambda^2  \\
b_{3,2}&=-(\tilde{D}_{2,1}\thd +2\tilde{\beta}\Omega) + (-2\tilde{\alpha}D_{3,2}\xi_3\Omega -2\xi_2\tilde{\beta}+\tilde{D}_{2,1}\xi_1)\lambda - \tilde{\alpha} \xi_3 (\tilde{D}_{2,1} \Omega+2 D_{3,2}\xi_2)\lambda^2 \\
b_{3,3}&=-\tilde{\alpha}( D_{3,1} \thd^2+D_{3,2} \Omega^2)+\tilde{\alpha} (2 D_{3,1}\xi_1 \thd -\tilde{D}_{2,1}\Omega^2-2 D_{3,2}\xi_2\Omega)\lambda\\
&-\tilde{\alpha} (\tilde{D}_{2,1}\xi_2\Omega+D_{3,2}\xi_2^2+D_{3,1}\xi_1^2)\lambda^2
\end{aligned}
\]
\section*{Appendix C}
In this section the numerical values used in the experiments are reported. All the experiments have the following common features 
\begin{itemize}
    \item $D_1=0.0411$, $D_3=0.8527$ (see, for instance, \cite{10.1063/1.4914122}),
    \item $D_2:=D_1+6.51 \lambda $ (perturbative setting),
    \item $\gamma=\sqrt{1-\Omega^2}$.
\end{itemize}
The remaining parameters have been chosen as follows
\begin{table}[h]
\begin{center}
\begin{tabular}{|c|c|c|c|c|c|c|}
\hline
Fig. & $D_2$ & $\tilde{\alpha}$ & $\tilde{\beta}_e$ & $\Omega$ & $\tilde{h}_2$ & $K$\\
\hline
	 2 & 0.1127 & 2 & 3 & 0.0676 & -0.4400 & 0.0700\\
\hline
	3 & 0.1127 & 2 & 3 & 0.0676 & -0.4400 & 0.0800\\
\hline
	4 & 0.0802 & 4 & $\times$ &  0.1799 & -1.1708 & 0.0524\\
\hline
	5 & 0.0802 & 4 & 3 & 0.1036 & -0.6741 & 0.0308\\
\hline
	6 & 0.0802 & 2 & 3 & 0.1036 & -0.6741 & 0.0308\\
\hline
	7 & 0.0802 & 2 & 3 & 0.1036 & -0.6741 & 0.0308\\
\hline
\end{tabular}
\label{tab:two}
\end{center}
\end{table}

\subsection*{Acknowledgements} This work has been supported by the Italian Ministry of University and Research, PRIN2020 funding program, grant number 2020PY8KTC. \\ 
The numerical simulations and the corresponding plots have been performed with GNU Octave \cite{oct}, whilst \cite{maxima} has been used for the algebraic manipulations. 

\bibliographystyle{alpha}
\bibliography{cqls.bib}

\end{document}